\documentclass{article}

\begin{document}

\title{Fredholm topology and enumerative
geometry: reflections on some words of Michael Atiyah}
\author{Simon Donaldson}
\maketitle













\newcommand{\bP}{{\bf P}}
\newcommand{\bC}{{\bf C}}
\newcommand{\bR}{{\bf R}}
\newcommand{\bQ}{{\bf Q}}
\newcommand{\bZ}{{\bf Z}}
\newcommand{\cO}{{\mathcal O}}
\newcommand{\cX}{{\mathcal X}}
\newcommand{\cU}{{\mathcal U}}
\newcommand{\cE}{{\mathcal E}}
\newcommand{\cB}{{\mathcal B}}
\newcommand{\cM}{{\mathcal M}}
\newcommand{\cL}{{\mathcal L}}
\newcommand{\oM}{{\overline{M}}}
\newcommand{\thalf}{{\textstyle \frac{1}{2}}}
\newcommand{\Cone}{{\rm Cone}}
\newcommand{\tih}{\tilde{h}}
\newcommand{\db}{\overline{\partial}}
\newcommand{\ocM}{\overline{{\mathcal M}}}
\newcommand{\ocN}{\overline{{\mathcal N}}}
\newcommand{\uGamma}{\underline{\Gamma}}
\newcommand{\maps}{{\rm Maps}}
\newcommand{\Diff}{{\rm Diff}}
\newcommand{\up}{\underline{p}}
\newcommand{\cW}{{\mathcal W}}
\newtheorem{premise}{Premise}
\section{Introduction}

{\it \lq\lq What we have here is a nonlinear way to use the index theorem''}.

\

It was my great good fortune to have many mathematical discussions with Michael
Atiyah, 
and I learnt an immense amount from him. Much of what he said to me I only
remember in general shape but in some cases I have more vivid memories and
the  quote above, uttered by Atiyah around 1982,  is such an example. (There
is another example at the beginning of Section 4 below). The first part of
this article is a general discussion motivated by this quote, and survey
of some developments over the past four decades. In the second part I will
outline  some new strands in these developments,  involving moduli spaces
of self-dual and complex structures. In addition to these brief reminiscences
the influence of Atiyah's ideas, work and general approach to mathematics
runs through the whole article.

\

\

{\it Index theory} can be seen as an extension of linear algebra to infinite
dimensions. Thus if $T:H_{1}\rightarrow H_{2}$ is a linear map between finite
dimensional vector spaces we have
   $$      {\rm dim}\  {\rm ker}\  T- {\rm dim}\ {\rm coker}\  T= {\rm dim}
\ H_{1}- {\rm dim}\  H_{2}. $$
If $H_{1}, H_{2}$ are infinite dimensional Banach spaces the right hand side
is not defined
but if $T$ is a {\it Fredholm} operator the kernel and cokernel are finite
dimensional (by definition) and the left hand side is defined: the index
of $T$. It is often useful to think, formally, of the index as the regularised
version of the difference of the two infinite dimensions. In particular it
is a deformation invariant on the space of Fredholm operators. This general
theory is not difficult: the significance of the {\it Atiyah-Singer Index
Theorem} is to connect with  geometry; finding the index for elliptic differential
operators over manifolds. 

{\it Differential topology}  is based on linear algebra. For example, if
we have submanifolds $P^{p}, Q^{q}$ of a manifold $M^{n}$ in general position
their intersection is a manifold of dimension $n-(p+q)$, just as happens
in the  linear case. Putting the two strands together, yields  a \lq\lq Fredholm
differential topology'' involving infinite dimensional manifolds and nonlinear
maps between them  with Fredholm derivatives. In the situations we consider
the general picture is that we have a  Banach manifold $\cB$, a vector bundle
with Banach space fibres $\cE\rightarrow \cB$ and a section
$\sigma$ which is locally represented by a nonlinear map with Fredholm derivative.
(That is, around a point $b\in \cB$ we choose a local coordinate chart and
a local trivialisation of the bundle $\cE$ so a section is represent by a
map from this chart to the fibre $\cE_{b}$ and we require that this be a
smooth map with Fredholm derivative.) 
The index of the derivative is independent of the local representation and
is formally the difference of the dimension of $\cB$ and the fibre dimension
of $\cE$. If the section $\sigma$ is transverse to the zero section then
the zero set $Z_{0}=\sigma^{-1}(0)$ is a manifold of dimension equal to the
index. The crucial issue is  the compactness of  $Z_{0}$. This is far from
automatic
since the ambient space $\cB$ is not compact  and, as we will see will typically
fail in a straightforward sense in cases of interest. But if for the moment
we assume the relevant compactness, then $Z_{0}$ carries a fundamental homology
class $\zeta=[Z_{0}]\in H_{\mu}(\cB)$ where $\mu$ is the index and this class
is a deformation invariant, independent of the choice of transverse section
$\sigma$ (for deformations through Fredholm sections, preserving compactness).
We assume that $\zeta$ can be defined in homology with rational coefficients
although that requires a discussion  of orientations. The proof of deformation
invariance is the same as in finite dimensions:  the standard differential
topological construction of the Poincar\'e dual of the Euler class of a vector
bundle.  So the upshot is that under suitable hypotheses there is a way to
define what is formally the (homology) Euler class of the infinite dimensional
bundle $\cE\rightarrow \cB$. 

In some cases of interest the homology groups of the infinite dimensional
ambient
space $\cB$ can be computed. In other cases one at least knows certain cohomology
classes which can be paired with the homology class $\zeta$ to produce
numerical invariants. That is, one has a  graded ring $R$ with a homomorphism
$R\rightarrow H^{*}(\cB)$ and the pairing gives an element of degree $\mu$
in $\check{R}= {\rm Hom}_{\bQ}(R, \bQ)$. 

These ideas, in the abstract, can be traced back a long way, certainly to
the 1965 paper \cite{kn:Smale} of Smale and are related to the older  Leray-Schauder
degree theory. The developments which are our focus here, starting around
1980, involve the application of these ideas to nonlinear differential equations
arising in geometry. Most of these developments fall into two broad lines.

\begin{enumerate}
\item Equations involving gauge fields, particularly over $4$-manifolds.
These include the Yang-Mills instanton equations (which were the context
for Atiyah's remark above) and the Seiberg-Witten equations. In the latter
case the theory has been extended to more sophisticated topological constructions
such as the  Bauer-Furuta invariants in stable homotopy \cite{kn:BF}. There
are other equations such as the Vafa-Witten equations which fit into the
Fredholm framework
but where compactness is only partially understood (with recent developments
in work of Taubes \cite{kn:Taubes2}). Similarly for various  analogues  over
manifolds of higher dimension, as discussed in \cite{kn:DT}: the main case
where compactness difficulties  are resolved occurs in the \lq\lq DT invariants''
for Calabi-Yau 3-folds and at present the resolution has to pass through
algebraic geometry \cite{kn:Thomas}. 

\item Pseudoholomorphic curves in symplectic manifolds. This goes back to
Gromov's 1986 paper and has become an enormous field. There are again some
analogous  equations for calibrated submanifolds in higher dimensions where
compactness is not adequately understood. 

\end{enumerate}

The situations in which the ideas sketched above can be usefully applied
(even making optimistic conjectures regarding compactness) are comparatively
rare and linked to \lq\lq low dimensional phenomena''.  There are many nonlinear
equations in differential geometry with Fredholm linearisations--very often
the linearisation is a variant of the Laplace operator. For one example we
can take closed geodesics $\gamma:S^{1}\rightarrow M$. But in most such cases
there is no hope of achieving compactness of the solution set: there could
be arbitrarily long  geodesics. On the other hand there are many situations
where some form of compactness holds but which involve {\it overdetermined}
equations.

To illustrate this, consider deformations of a compact complex submanifold
$P$ in a fixed complex manifold $M$.  Let $N_{P} = TM\vert_{P}/TP$ be the
normal bundle-- a holomorphic bundle over $P$.
The linearised equation is given by the $\db$-operator on sections of $N$
and the solutions, which give the infinitesimal deformations of $P$, are
 the holomorphic sections $H^{0}(P;N_{P})$. But it may not be possible to
extend these to genuine deformations; there are potential obstructions in
$H^{1}(P;N)$.  In the case of curves, when ${\rm dim}_{\bC} P=1$, this fits
perfectly into our Fredholm picture (as in item (2) above). The cohomology
groups $H^{0}(P,N_{P}), H^{1}(P,N_{P})$ are the kernel and cokernel of the
linearised operator and the (real) index is twice the Euler characteristic
$\mu= 2( {\rm dim}\  H^{0}- {\rm dim}\  H^{1})$. If we are in a transverse
situation
then  $H^{1}$ vanishes and the space $Z_{0}$ of holomorphic curves near to
$P$ is a (real) manifold of dimension $\mu$ (and in fact a complex manifold).
But it might happen that we are not in a transverse situation and $Z_{0}$
could
have  some very different structure: it could be singular or a manifold of
dimension greater than $\mu$. Then we can follow two paths, one differential
geometric and one algebro-geometric. (One expects these to reach the same
endpoint, although technically this may be highly non-trivial, in general.)

\begin{itemize}
\item Differential geometrically, we can perturb our equations
in some way, for example by perturbing the complex structure on $X$ to an
almost-complex structure, so that the perturbed equation is transverse and
we get a solution set $Z$ of the perturbed equation which is a manifold of
dimension $\mu$.
\item Algebro-geometrically, one can add additional structure
to the  original  set $Z_{0}$ which enables one to define a {\it virtual
fundamental class} in $H_{\mu}(Z_{0})$. In the simplest case, when the potential
obstructions do not actually occur and $Z_{0}$ is a manifold of real dimension
$2 {\rm dim} H^{0}(N)$, one considers  the vector bundle over $Z_{0}$ formed
by the cohomology groups $H^{1}(P;N_{P})$. The virtual fundamental class
is the Poincar\'e dual of the Euler class of this bundle, which lies in $H_{\mu}(Z_{0})$.

\end{itemize}

 The picture is fundamentally different when $P$ has dimension greater than
$1$. From the differential geometric point of view, the equations defining
a complex submanifold become overdetermined and if we perturb the complex
structure on $M$ to a generic almost-complex structure we expect there to
be no solutions. From an algebro-geometric point view we still have a deformation
theory,  with infinitesimal deformations in $H^{0}(P, N_{P})$ and potential
obstructions in $H^{1}(P, N_{P})$.  There is a \lq\lq Kuranishi model'' of
 a neighbourhood in the space  of submanifolds as the solutions of
$h_{1}$ equations in $h_{0}$ complex variables where $h_{i}={\rm dim}\  H^{i}$,
so it might seem reasonable to think of the \lq\lq expected'' complex dimension
of this neighbourhood as $h_{0}-h_{1}$, as before. The difference is that
there are now higher cohomology groups $H^{i}(P; N_{P})$ for $i\geq 2$.
The Euler characteristic $\sum (-1)^{i} h_{i}$ is a deformation invariant
but, without some control of the higher cohomology, the difference $h_{0}-h_{1}$
is not and the expected dimension could be different at different points
of $Z_{0}$.

For another example, consider the space $V$ of conjugacy classes of irreducible
representations in $SU(r)$ of the fundamental group $\pi=\pi_{1}(M)$ of a
compact oriented manifold $M$. This set $V$ has the structure of a real algebraic
variety and clearly only depends on the group $\pi$. When $M$ has dimension
$3$ there is a sense in which the \lq\lq expected dimension'' of $V$ is $0$,
even though the actual dimension of $V$ could be very different. This is
the idea behind the Casson invariant (for homology $3$-spheres) which \lq\lq
counts'' the points in $V$. To fit this into our general framework (as was
done by Taubes \cite{kn:Taubes1}), we consider $V$ as the set of isomorphism
irreducible
flat connections on an $SU(r)$ bundle $E\rightarrow M$. The infinite dimensional
space $\cB$ is the quotient of the space of all connections by the gauge
group ${\rm Aut}(E)$. The curvature $F(A)$ of a connection $A$ can be viewed
as a cotangent vector in the space of connections and we get a Fredholm section
$\sigma$ of the cotangent bundle $T^{*}\cB\rightarrow \cB$ whose zero set
is identified with $V$. Then the Casson invariant is one half the homology
class we discussed above: i.e. half the number of zeros (counted with suitable
signs) of a generic perturbation of $\sigma$.  From another point of view,
the deformation theory of a flat connection $A$ over a manifold $M$ can be
discussed through the de Rham complex:
\begin{equation} \Omega^{0}( {\rm ad}\  E)\stackrel{d_{A}}{\rightarrow}\Omega^{1}({\rm
ad}\  E)
\stackrel{d_{A}}{\rightarrow}   \Omega^{2}( {\rm ad}\  E)\dots  \end{equation}
Here ${\rm ad}\  E$ is the bundle of Lie algebras associated to the adjoint
representation and $d_{A}$ is the coupled exterior derivative defined by
$A$. The cohomology groups $H^{i}_{A}$ of the complex are the cohomology
groups of the manifold $M$ in the local coefficient system ${\rm ad} E$.
If, as we are assuming, the connection is irreducible the group $H^{0}_{A}$
vanishes. The group $H^{1}_{A}$ corresponds to infinitesimal deformations
of $A$ and there are potential obstructions in $H^{2}_{A}$. Poincar\'e duality
gives an dual pairing between $H^{i}_{A}$ and  $H^{n-i}_{A}$ so when $M$
is a $3$-manifold there are just two non-zero cohomology groups $H^{1}_{A},
H^{2}_{A}$ which are dual,  and this fits in with the fact that the expected
dimension of $V$ is zero. For a higher dimensional manifold we are in the
same position as in the previous example: there are higher dimensional cohomology
groups $H^{3}_{A}, \dots$ and we don't have a way to define an expected dimension
of $V$.

There is a very similar discussion for holomorphic vector bundles over compact
complex manifolds. The deformation theory of a holomorphic bundle $E\rightarrow
X$ (with fixed determinant) is
discussed through the cohomology groups $H^{i}(X, {\rm End}_{0} E))$, where
${\rm End}_{0} E $ is the bundle of trace-free endomorphisms. The group $H^{0}$
is related to the \lq\lq stability'' of the bundle and we can assume it is
zero for this discussion.  The groups $H^{1}, H^{2}$
give infinitesimal deformations and obstructions, as before. To have a virtual
dimension we need control of the higher cohomology groups $H^{i}$ which would
typically come through vanishing: $H^{i}=0$ for $i\geq 3$. This is automatic
if $X$ is a complex surface (and then the holomorphic bundles are intimately
related to Yang-Mills instantons). If $X$ is a Calabi-Yau $3$-fold then Serre
duality gives the vanishing of $H^{3}$ and the picture is closely analogous
to the previous discussion of the Casson invariant: this is the starting
point for the definition of \lq\lq DT invariants'' counting holomorphic bundles
(or more generally sheaves) over $X$ \cite{kn:Thomas}.  More generally, we
get vanishing of $H^{3}$ for any $3$-fold $X$ with a nonzero section of $K_{X}^{-1}$.
So, for example, moduli spaces of bundles over $X=\bC\bP^{3}$ have an \lq\lq
expected dimension (which is also related to Yang-Mills instantons on $S^{4}$,
through twistor theory). But there is no obvious way to define an expected
dimension for moduli spaces of bundles over $\bC\bP^{n}$ for $n\geq 4$.

These examples illustrate the special and low-dimensional nature of these
Fredholm and virtual fundamental class techniques. There are some situations
in which these ideas can be applied  in less direct ways. 
\begin{itemize}
\item  Furuta and Ohta define a Casson-like invariant through representations
of the fundamental group of $4$-manifolds by interpreting  flat connections
as Yang-Mills instantons \cite{kn:FO}. 
\item  Borisov and Joyce \cite{kn:BJ} and Cao and Leung \cite{kn:Cao}  define
enumerative invariants for sheaves over Calab-Yau 4-folds by, from a differential
geometer's point of view,  interpreting these
within ${\rm Spin}(7)$-geometry.  
\end{itemize}

\section{Towards enumerative theories for structures on manifolds}

 In this section we discuss two other cases where these Fredholm/virtual
fundamental
class techniques might be applicable; different in nature from those considered
above. We begin with some general background in Subsection 2.1 and then review
the geometric setting for these two examples in Subsection  2.2. 

\subsection{Miller-Mumford-Morita classes}

 A general picture is that we want to consider some kind of structure
on a compact oriented $n$-manifold $M$, which we can take to be given by
a tensor field $\tau$ (for example a Riemannian metric). So we have an infinite
dimensional space $\cX$ of these tensors which is acted on by the group $\Diff$
of orientation preserving diffeomorphisms of $M$ and we form the quotient
space $\cB=\cX/\Diff$.  Suppose initially that the action is free,  so that
$\cB$ is an infinite dimensional manifold. The space $\cX$ is a principal
$\Diff$ bundle over $\cB$ and we have a classifying map
$$    g: \cB\rightarrow B\Diff, $$
and  $ g^{*}: H^{*}(B\Diff)\rightarrow H^{*}(\cB)$. In examples the assumption
that $\Diff$ acts freely is too restrictive and we should weaken it to the
assumption that stabilisers are finite. Then $\cB$ is an infinite dimensional
orbifold and if we use rational coefficents we still have a map $g^{*}$.

The cohomology of $B\Diff$ might not be accessible but there is a standard
way to produce classes in it. Let $\gamma$ be a $p+n$-dimensional rational
characteristic class of $SL(n;\bR)$ (or equivalently $SO(n)$) with $p\geq
n$. There is a universal bundle
$\cU\rightarrow B\Diff$ with fibre $M$ and a tangent bundle along the fibres
$T_{V}\rightarrow \cU$, so we  get $\gamma(T_{V}) \in H^{p+n}(\cU)$ and integration
over the fibre of $\cU\rightarrow B\Diff$ gives us a class in  $ H^{p}(B\Diff)$.
Finally we pull-back by $g$ to get 
$I(\gamma)\in  H^{p}(\cB) $.
We could define $I(\gamma)$ without going through $B\Diff$ and the universal
bundle by considering directly an orbifold bundle with fibre $M$ over $\cB$.

The classical example is when $n=2$, for Riemannian metrics on  an oriented
surface. The characteristic classes in question are just polynomials in the
Euler class, or first Chern class, $c_{1}\in H^{2}$ and
we get {\it  Miller-Mumford-Morita classes} $I(c_{1}^{p+1})\in H^{2p}(\cB)$.
 In other words, if we form a graded ring 
 $${\mathcal R}= \bQ[\sigma_{1}, \sigma_{2}, \dots]$$
  freely generated by objects $\sigma_{p}$ in degree $2p-2$ then we have
a homomorphism  
  $$   \Gamma:{\mathcal R}\rightarrow H^{*}(\cB),$$
  taking $\sigma_{p}$ to $I(c_{1}^{p+1})$. 
  This extends in the obvious way to define a ring ${\mathcal R}_{SO(n)}$
in each dimension
$n$ taking characteristic classes which are polynomials in the Pontryagin
classes and, for even $n$, the Euler class. There have been many recent advances
in geometric topology,  in the understanding of these generalised Miller-Mumford-Morita
(MMM) classes and $H^{*}(B\Diff)$, see the survey \cite{kn:GRW} for example.
(These classes are also called \lq\lq tautological classes'' in the literature.)

There are parallel constructions using $K$-theory. Suppose that $n$ is even,
the manifold $M$  is a spin manifold and we have a spin structure on the
vertical tangent bundle $T_{v}$. Take a Euclidean metric on $T_{v}$ so the
fibres of $\cU\rightarrow B\Diff$ become Riemannian spin manifolds. Let $\rho$
be a representation of ${\rm Spin}(n)$. We get an associated bundle over
the fibres and a coupled Dirac operator $D_{\rho}$. The index of the family
defines a class  ${\rm ind} D_{\rho}$ in $K(B\Diff)$ and the Atiyah-Singer
theory gives a formula for the Chern character ${\rm ch}\ ({\rm ind}\  D_{\rho})\in
H^{*}(B\Diff)$ in terms of the MMM classes. Conversely we can express all
the MMM classes in terms of the ${\rm ind} D_{\rho}$. (Note: For questions
involving rational cohomology the  spin assumption is not important because
the $D_{\rho}$ exist  for at least half of the representations $\rho$,  which
suffice to generate all the MMM classes. Also in the construction above we
should strictly restrict to compact subsets of $B\Diff$, but this also will
not be important for us.)

There is an extension of this discussion in the case when the tensors $\tau$
we are considering give a reduction of the structure group of $TM$ to some
subgroup $G\subset SL(n,\bR)$ which is not homotopy equivalent to $SO(n)$.
Then we can start with other characteristic classes of $G$ to get a ring
${\mathcal R}_{G}$ which maps to   $H^{*}(\cB)$. Another variant is when
we have as part of our structure
a  symplectic form $\omega$  on the manifold $M$ which we take as fixed and
reduce the symmetry group to the symplectic diffeomorphisms ${\rm SDiff }
(M, \omega) \subset \Diff$.

With this background in place we come to the main point. Suppose that we
have a suitable differential equation for structures $\tau$. That is, we
have some vector bundle $\cE\rightarrow \cB$, a Fredholm section $\sigma$
of index $\mu$, as considered
in the previous section, and our equation is $\sigma(\tau)=0$. Then, modulo
compactness and orientation questions, we can fit into that general framework
and define an element of degree $\mu$ in the dual $\check{{\mathcal R}}_{SO(n)}$
of the
ring ${\mathcal R}_{SO(n)}$ (or $\check{{\mathcal R}}_{G}$ for a structure
group $G$). 

\

\

In  the classical case, when $\tau$ is a Riemannian metric on a surface of
genus $g\geq 2$ we  take  $\sigma(\tau)= K_{\tau}+1$ where $K_{\tau}$ is
the Gauss curvature. So the solutions of our equation are metrics of curvature
$-1$, the index is $6g-6$ and the zero set $Z\subset \cB$ is the moduli space
$\cM_{g}$. By the uniformisation theorem we achieve the same end  by taking
$\tau$ to be an almost-complex structure, with no equation, and in that
set-up we just have $\cB= \cM_{g}$, so $\cB$ is finite-dimensional. Since
$\cX$ is contractible the map $g^{*}$ is an isomorphism on rational cohomology.
There is a huge literature on
$H^{*}(\cM_{g})$, the MMM classes and the integrals of MMM classes over $\cM_{g}$.
On the $K$-theory side, the irreducible representations of $Spin(2)=S^{1}$
are labelled by an integer $r$. Taking for simplicity $r=1-2q$ for integers
$q\geq 1$, the corresponding virtual bundle over $\cM_{g}$ is the vector
bundle which assigns to each complex structure the vector space $H^{0}(M,K_{M}^{q})$.

\subsection{Two structures in real dimension $4$}

Before going further let us emphasise that, in the context of possible enumerative
theories, the material  we are discussing below is largely speculative and
on a very different footing from the well-established theories we reviewed
in Section 1.  
\begin{itemize}\item
Our first example is {\it self-dual conformal structures}, for which  a basic
reference is the paper \cite{kn:AHS} of Atiyah, Hitchin and Singer. Recall
that if $(M,g)$ is an oriented Riemannian $4$-manifold the Weyl tensor decomposes
as $W=W_{+}\oplus W_{-}$ where $W_{\pm}\in s^{2}_{0}(\Lambda^{2}_{\pm})$,
 the trace-free symmetric $2$-tensors on the $3$-dimensional $\Lambda^{2}_{\pm}$.
The Weyl curvature is conformally invariant and a conformal structure with
$W_{-}=0$ is called self-dual. In our set-up we take $\cX$ to be the space
of conformal structures $\tau$ and the section to be that induced by $W_{-}$,
so the zero set $Z_{0}$ is the moduli space of self-dual conformal structures
modulo diffeomorphism. This is a Fredholm section of index 
\begin{equation}  {\rm index}= \frac{1}{2} (29 \sigma(M)- 15 e(M)), \end{equation}
where $\sigma(M), e(M)$ are  the signature and  Euler characteristic.
The  deformation theory is worked out in \cite{kn:KK}. There is an elliptic
complex
\begin{equation}  \Gamma(TM)\rightarrow \Gamma(s^{2}_{0}T^{*}M)\rightarrow
\Gamma(s^{2}_{0}\Lambda^{2}_{-})\end{equation}
where the first term consists of vector fields on $M$ (the Lie algebra of
$\Diff$), the second consists of deformations of a conformal structure (the
tangent space of $\cX$) and the third is the space where $W_{-}$ lives {\it
a priori} (the fibre of the bundle $\cE$).

The MMM classes in this situation lie in dimensions $0$ modulo $4$ so  for
a manifold $M$ with $e=3\sigma\ \ \ {\rm mod}\  8$ one might hope to define
 pairings with a virtual fundamental class.  

\item 
The second example is {\it complex surfaces}. Ignoring for the moment some
fundamental difficulties, we consider a compact oriented $4$-manifold $M$
and the space $\cX$ of almost-complex structures $\tau$. The section $\sigma$
of a bundle over $\cB$ is induced by the Nijenhuis tensor
$N(\tau)$ of an almost complex structure, which lies in $\Lambda^{0,2}(T)$
and this is a Fredholm section of (real) index 
\begin {equation} \mu= 2\chi(TM) = 2 (10\chi(S)- 2c_{1}^{2}(S)). \end{equation}
Here we are writing $\chi(S)$ for the holomorphic Euler characteristic which
is equal, by the Riemann-Roch formula, to $(c_{1}^{2}+c_{2})/12$. Using standard
formulae we can also write $\mu=-(3e+7\sigma)$.

As a set, the zeros of our section correspond to equivalence classes of complex
structures on $M$, which we would like to call the \lq\lq moduli space''
of complex structures. The fundamental difficulty we encounter is the well-known
one that this is not in general a good object: the natural topology on the
set may not be Hausdorff. We will come back to this in the next section.

For a general complex manifold $X$  the Kodaira-Spencer-Kuranishi theory
describes the \lq\lq versal deformation'' in terms of the cohomology groups
 $H^{p}(TX)$. As in the examples from the previous section, the special feature
of complex dimensions $1,2$ are reflected in the fact that there  are no
higher cohomology groups $H^{i}$ for $i\geq 3$. 

In this situation the  structure group is $GL(2,\bC)$ and the characteristic
classes are generated by $c_{1}, c_{2}$. The relevant ring is
$$ {\mathcal R}_{U(2)}= \bQ[ \sigma_{p q}], $$
with generators $\sigma_{pq}$ of degree $2(p+2q-2)$ corresponding to $c_{1}^{p}
c_{2}^{q}$. One might--- in suitable situations---hope to define pairings
of ${\mathcal R}_{U(2)}$ with a virtual fundamental class. 
( Even for moderate values of the index $\mu$ this would give a  large  collection
of
numbers.  
For example, the dimension of the degree $8$ part of ${\mathcal R}_{U(2)}$
is $30$.)
\end{itemize}
\

\

It has to be said that there is not the same clearcut motivation for developing
such enumerative theories as there is in other situations,  with the construction
of 4-manifold and symplectic topology invariants. It is also not so easy
to find natural deformations of the equations (although we will encounter
something on those lines in the next section). From the point of view
of this article we could say that the motivation is that, as we have explained,
the situations where this kind of nonlinear Fredholm theory is possibly relevant
are comparatively
rare and special, hence precious,  and  one wants to understand them as
far as possible.

  \section{Compactification}

  \

  {\it \lq\lq You need to be careful compactifying moduli spaces: people
spend their lives doing that''}.

\

This is our second quote from Michael Atiyah. The context was that we were
discussing the \lq\lq Uhlenbeck compactification'' of moduli spaces of instantons.
For the immediate purposes then, given the fundamental analytical results
of Uhlenbeck, this was quite straightforward to define. But Atiyah's remark
holds true in the sense that understanding in detail the structure of the
Uhlenbeck compactification is crucial in establishing deep properties of
the instanton invariants of 4-manifolds such as  Witten's conjecture on the
relation with Seiberg-Witten invariants, as in the work of Feehan and Leness
\cite{kn:FL}, and this is something which is still not fully understood.

\

In the context of this article Atiyah's remark points to the core of the
matter.  The zero set $Z_{0}$ of our Fredholm section will usually not be
compact
so does not carry a homology class and the evaluation of cohomology classes
has no clear meaning.   What we would like to do is to compactify the zero
set in such a way that the relevant cohomology classes and the deformation
theory extend over the compactification.

\

There is not much known  about compactification of moduli spaces of self-dual
conformal structures. If we have a sequence of such structures $\tau_{i}$
on a fixed smooth 4-manifold $M$ we would like to identify some kind of geometric
limit
of a subsequence. If we suppose that within each conformal class there are
metrics $g_{i}$ with constant scalar curvature and bounded
Sobolev constant then there are results of Tian and Viaclovsky \cite{kn:TV}.
But one can say that the construction of a rigorous general theory seems,
at best, very far off. There are some  explicit examples known of connected
components of moduli
spaces and we mention two.

\

{\bf Example 1}

\

  For the manifold $M=\bC\bP^{2}\sharp \bC\bP^{2}$ there is a component
of the moduli space constructed by Poon \cite{kn:Poon} which is an open interval
$(0,1)$. (The index here is $-1$ but the structures have 2-torus symmetry
which modifies the \lq\lq expected dimension'' discussion, so in this context
the dimension is the expected one.) The natural compactification is the closed
interval $\overline{Z}=[0,1]$. The geometric meaning is that the point $0\in
\overline{Z}$ can be thought of as corresponding to a wedge $\bC\bP^{2}\vee\bC\bP^{2}$
of two copies of $\bC\bP^{2}$ with its Fubini-Study metric, joined at a point.
The convergence to this singular limit is realised by a sequence of conformal
structures on the connected sum with the \lq\lq neck'' shrinking to zero
size. The other end point $1$ is similar. There is a well-known  Eguchi-Hanson
metric on the tangent bundle of the $2$-sphere, which is self-dual and asymptotically
locally Euclidean. The conformal 1-point compactification of this gives a
compact self dual orbifold $T$ with one singular point and the end-point
$1$ corresponds to the wedge $T\vee T$ of two copies glued at their singular
points. The convergence is realised by shrinking the neck of a \lq\lq connected
sum'' of these orbifolds. (See the discussion in \cite{kn:DF}.)

\

{\bf Example 2}

\

 Let $M$ be the $K3$ manifold with the opposite of its standard orientation.
The Calabi-Yau metrics on $M$ define self-dual conformal structures and the
Torelli theorem for K3 surfaces shows that a connected component of the moduli
space has an explicit description   $U/\Gamma$. Here  $U$ is a certain dense
open
set in the negative Grassmannian $Gr^{-}(19,3)$ of negative $3$-dimensional
subspaces in the indefinite space $\bR^{19,3}$  and $\Gamma$ is a subgroup
of the isometry group of the $K3$ lattice  $\Lambda_{K3}\subset \bR^{3,19}$.
The complement of $U$ in the Grassmannian is a union of explicit  codimension-$3$
sets and these correspond to structures with mild orbifold singularities.
As a first step towards a compactification we can add these points to get
a moduli space 
$Z_{0}= Gr^{-}(19,3)/\Gamma$.
  
The index formula (2) gives a virtual dimension 52 whereas the dimension
of $Gr^{-}(19,3)$ is $3.19=57$. The difference is accounted for by the fact
that  the linearised operator has a  cokernel: i.e. the cohomology  $H^{2}$
of the deformation complex (3) is nonzero. The Calabi-Yau metric induces
a flat connection on the bundle $s^{2}_{0}\Lambda^{2}_{-}$ and the space
$H^{2}$ can be identified with the $5$-dimensional space of parallel sections.
This gives an explicit description of the obstruction bundle over $Z_{0}$
as the quotient by $\Gamma$ of the rank $5$ bundle $s^{2}_{0}(V)$, where
$V$ is the tautological $\bR^{3}$ bundle over the Grassmannian. So a generic
perturbation of the moduli space  corresponds to the zero set $Z\subset Z_{0}$
of a generic section of $s^{2}_{0}(V)$.
 To develop this further we would have to consider a suitable compactification
of the moduli space $Z_{0}$, which is an important topic for other purposes,
but there are some other interesting points that arise. 
 
 Recall that the indefinite orthogonal group $O(19,3)$ has four connected
components, corresponding to the action on the orientations of positive and
negative subspaces. So we have two homomorphisms $o_{+},o_{-}: O(19,3)\rightarrow
\{\pm 1\}$. It is known that the group $\Gamma$ is the intersection of $O(\Lambda_{K3})$
with the kernel of $o_{-}$: the diffeomorphisms of $M$ preserve the orientation
of the negative subspace and this is a fundamental phenomenon in smooth
4-manifold theory \cite{kn:Donaldson1}. On the other hand the homomorphism
$o_{+}$ is non-trivial on $\Gamma$: the reflections by classes of square
$2$ are realised by generalised Dehn twists. These observations show that
the bundle $s^{2}_{0}(V)$ over $Z_{0}$ is orientable but the space $Z_{0}$
is not, so we cannot define a virtual fundamental class in rational homology
even ignoring compactness. But we will suggest a variant of the set-up which
gets around this.

In both the $SU(2)$-instanton and Seiberg-Witten theories the orientation
of moduli spaces is governed by the orientation of the positive and negative
parts of the second cohomology of the underlying 4-manifold. Consider pairs
$(\tau,A)$ consisting of a self-dual conformal structure and an anti-self-dual
connection $A$ of Chern class $k$. The discussion above shows that, at least
in this case of the K3 manifold, the moduli space of these pairs {\it is}
orientable. It  has virtual dimension $8k-60+52=8k-8$ (which is divisible
by $4$). One can also  also deform to an interesting coupled system. We assume
for simplicity that that we have chosen some way to fix a metric within each
conformal class and consider the equations for pairs $(\tau,A)$:
 $$   F^{+}_{A}=0  \ \ \ \,\ \ \  W^{-}(\tau)= \epsilon F^{-}_{A}* F^{-}_{A},
$$
where $\epsilon$ is a real parameter and  $*$ denotes the quadratic map combining
the Killing form on the Lie algebra with $\Lambda^{2}_{-}\otimes \Lambda^{2}_{-}\rightarrow
s^{2}_{0}(\Lambda^{2}_{-})$. It seems likely that for generic $\epsilon$
the moduli space of solutions of this coupled system is an orbifold of dimension
$8k-8$.  There is a similar discussion for the Seiberg-Witten case. The constructions
are related to the Seiberg-Witten invariants for families  \cite{kn:Nak},
\cite{kn:LiLiu}. It would be interesting to pursue a general study of these
orientation questions for self-dual structures.

 Another obvious issue in this K3 discussion is that the rank of the obstruction
bundle $s^{2}_{0}V$ is odd and in a standard situation  the Euler class of
a bundle of odd rank vanishes in rational cohomology. However it is possible
that the right  treatment of the compactification could allow non-trivial
pairings.

\section{Surfaces of general type}

We now turn to the main technical topic of this article, considering complex
structures on 4-manifolds defining  surfaces of \lq\lq general type'', which
are the analogues of complex curves of genus two or more. There is a huge
literature about these and in particular there is a well-developed Koll\'ar,
Shepherd-Barron, Alexeev (KSBA) theory of moduli space compactification.

 To set up the basic picture differential geometrically we can start with
a compact symplectic $4$-manifold $(X,\omega)$. The symplectic structure
defines a first Chern class $c_{1}\in H^{2}(X)$ and we assume that $c_{1}=-[\omega]$.
Consider the space
${\cal X}$ of almost-complex structures  on $X$, algebraically compatible
with $\omega$, and let $\cB$ be the quotient by the symplectic diffeomorphism
group. A compatible almost-complex structure defines a Riemannian metric
$g(J,\omega)$ on $X$ and this means that there is no difficulty in forming
the quotient space. One can define a \lq\lq Hermitian scalar curvature''
$S(J,\omega)$ which reduces to the ordinary scalar curvature when $J$ is
integrable and which has the property that
$$ \int_{X} \left(S(J,\omega)+ 4\right)  \omega^{2}= 0. $$

We consider the section $\sigma$ of a bundle over $\cX$ which corresponds
to the pair of tensors $(N(J), S(J,\omega)+4)$,  where $N(J)$ is the Nijenhuis
tensor.
 Thus a zero of $\sigma$ gives a K\"ahler structure on $X$ in the class $-c_{1}$
 with scalar curvature $-4$ and a standard integral identity then shows that
this metric is in fact a {\it K\"ahler-Einstein} metric, with Ricci curvature
$-g(J,\omega)$. Thus the space $Z_{0}\subset \cB$ is the moduli space of
K\"ahler-Einstein structures. This is a Fredholm section of index $\mu$ where
$\mu$ is given by the same formula (4). If, for simplicity, we assume that
$H^{1}(X,\bR)=0$ the relevant deformation complex is

\begin{equation}  C^{\infty}_{0} \rightarrow \Gamma( s_{\bC}^{2}T)\rightarrow
\Gamma(\Lambda^{0,2}\otimes T) \oplus C^{\infty}_{0}, \end{equation}
where $C^{\infty}_{0}$ denotes real valued functions of integral zero and
$s_{\bC}^{2}T$ is the symmetric square of the tangent bundle, regarded as
a complex vector bundle using $J$. 

A complex surface which admits a \lq\lq negative'' K\"ahler-Einstein metric,
 ${\rm Ricci}=-g$,  is of general type and the converse is almost true. Let
$Y$ be a smooth complex surface of general type. It might be that the canonical
bundle $K_{Y}$ is ample and in that case $Y$ admits a unique K\"ahler-Einstein
metric by the theorem of Aubin and Yau so we can take $X=Y$ above. But in
the theory it is best to consider moduli spaces of all structures with fixed
numerical invariants
$(c_{1}^{2}, \chi)$, so there could be {\it different} underlying symplectic
4-manifolds $(X,\omega)$ (although no example is known). If $K_{Y}$ is not
ample there is a {\it canonical model} obtained by contracting all $-1$ and
$-2$ curves.
This is an orbifold with ADE singularities and carries a unique orbifold
K\"ahler-Einstein metric by an extension of the Aubin-Yau theorem \cite{kn:Kobayashi}.
 In short, before  compactification,  we should strictly run the differential
geometric discussion above with possibly a finite number of different symplectic
$4$-manifolds or orbifolds,  but we will not go into this further.

The KSBA theory produces a compactified moduli space  $\ocM_{a,b}$  of  surfaces
of general type with $c_{1}^{2}=a, \chi=b$. (The original references are
\cite{kn:Alexeev}, \cite{kn:KSB} and there is a helpful  exposition in \cite{kn:Hacking}.)
It is analogous to the Deligne-Mumford compactification $\ocM_{g}$ of curves
of genus $g\geq 2$. Berman and Guenancia show that the singular  surfaces
represented in the compactification are precisely those which admit K\"ahler-Einstein
metrics, with a suitable technical definition of what that means in the singular
case \cite{kn:Berman}. This is analogous to the hyperbolic
geometry description of the Deligne-Mumford space, but the metrics are not
necessarily complete and the picture is much more complicated. There are
many interesting questions about the asymptotics of these metrics. In any
case, while in the future it may be possible to proceed in a more differential
geometrical fashion,
we will now switch to a purely algebro-geometric point of view. Our
discussion is based on the following foundational premise.

\

%
{\bf PREMISE}
{\it 

\begin{enumerate}
\item  The virtual fundamental class theory of Behrend and Fantechi  \cite{kn:BFan}
and Li and Tian \cite{kn:LT}, or some variant of that,
can be used to define a class $\zeta \in H_{\mu}(\ocM_{a,b},\bQ)$ where $\mu=\mu(a,b)=20b-4a$.
\item The MMM classes extend to $H^{*}(\ocM_{a,b},\bQ)$.
\end{enumerate} }

The author does not have the expert knowledge required to write a useful
  discussion of this premise and is certainly not suggesting that it must
be true: the main point here is to raise the questions. The author's impression
is that this statements should be true at least for moduli spaces in which
the surfaces involved are not too badly singular,  as in the example we study
in Section 5 below. 
In any event, assuming---for now---the premise, we immediately get the existence
of elements $\rho_{a,b}\in \check{{\mathcal R}}_{U(2)} $ for each $(a,b)$
such that $\mu(a,b)\geq 0$, which are the    main  point of this article.

\

{\bf  Remarks}

\begin{enumerate}
\item One does not have to go far to encounter cases where the virtual moduli
space dimension is different from the actual dimension.  For example, for
a smooth surface $S$ of degree $d\geq 5$ in $\bC\bP^{3}=\bP(U)$ standard
exact sequences show that $H^{2}(TS)$ is isomorphic to the cokernel of the
natural map
$$   s^{d-4}(U)\rightarrow s^{d-5}(U)\otimes U$$
which has dimension $(d-2)(d-3)(d-5)/2$ and is non-trivial if $d\geq 6$.
The actual dimension of the moduli space is $$\frac{(d+1)(d+2)(d+3)}{6}-16,
$$  while the virtual dimension is less by the dimension of $H^{2}(TS)$.
Catanese \cite{kn:Catanese} and Manetti \cite{kn:Manetti} have shown that
the moduli spaces $\cM_{a,b}$ can have a  large number of components of different
dimensions.  
\item There is a well-known \lq\lq geography'' of surfaces of general type,
with a region ${\cal S}$ in the $(a,b)$-plane outside which the moduli spaces
 $\cM_{a,b}$
are empty. (See \cite{kn:BPV}, VII.9, for example.)
The line $\mu(a,b)= 0$ cuts through the middle of this region ${\cal S}$.
For many interesting
surfaces $\mu(a,b)<0$ and the virtual class theory evaporates. From the differential
geometric point of view one can ask the question, when does a given symplectic
4-manifold admit a K\"ahler structure? One might try to develop some analytical
scheme to prove existence. But it is harder to imagine how  such a scheme
could find these surfaces with $\mu<0$,  since at the index level the solutions
\lq\lq ought not to exist''.

\item  The obstruction space for deformations of a compact complex manifold
$V$ of any dimension is always $H^{2}(TV)$. A special feature of complex
dimension $2$ is that this space has a direct geometric interpretation: the
Serre dual is then $H^{0}(T^{*}V\otimes K_{V})$ and elements of this give
singular holomorphic foliations of the surface $V$. This is analogous to
the description, in complex dimension $1$, of the cotangent bundle of $\cM_{g}$
in terms of quadratic differentials. 

\item In regard to to the second item in the \lq\lq Premise'' the two-dimensional
MMM classes exist in  $H^{2}(\ocM_{a,b};\bQ)$ and are studied in the literature.

 In terms of line bundles, there are Knudsen-Mumford line bundles $\cL_{0},
\cL_{2}\rightarrow \ocM_{a,b}$ such that:
\begin{equation}  {\rm det} \ \pi_{*}(K^{p})= \cL_{2}^{N_{2}(p)}\otimes \cL_{0}^{N_{0}(p)}\end{equation}
where, writing $\up=p-1/2$,  
\begin{equation}  N_{2}(p)= \frac{\up^{3}}{3}- \frac{\up}{12}  \ \ ,\ \ 
 N_{0}(p)=- 2\up \end{equation} 
 The left hand side in (6) refers to the relative canonical bundle of the
family
$\cU\rightarrow \ocM_{a,b}$. Note that for a general ample line bundle $L$
over a surface the Knudsen-Mumford theory gives four lines $\cL_{3},\cL_{2},\cL_{1},
\cL_{0}$ such that
$$    {\rm det}\ \pi_{*}(L^{p})= \cL_{3}^{n_{3}(p)}\cL_{2}^{n_{2}(p)}\cL_{1}^{n_{1}(p)}\cL_{0}^{n_{0}(p)},
$$
where $n_{i}(p)= \left(\begin{array}{cc} p\\i\end{array}\right)$, but for
the canonical bundle Serre duality implies that $\cL_{3}= \cL_{2}^{2},\cL_{1}=\cL_{0}^{-2}$
and we get the expressions (6),(7).

 To be more precise, these are all orbifold or $\bQ$-line bundles, due to
the possible presence of finite automorphism groups. The Grothendieck-Riemann-Roch
theorem shows that
\begin{equation} c_{1}(\cL_{2})= -\frac{1}{2} I(c_{1}^{3})\ \ \ ,\ \ \ \
 c_{1}(\cL_{0})= 24 I(c_{1}c_{2})\end{equation}
The line bundle $\cL_{3}=\cL^{2}_{2}$ is known as the CM line bundle. Patakfalvi
and Xu \cite{kn:PXu} show that it is an ample line bundle over $\ocM_{a,b}$,
and it follows the MMM class  $I(c_{1}^{3})$ is non-zero in $H^{2}(\ocM_{a,b};\bQ)$.

In the opposite direction, Randall-Williams shows that on a moduli space
of {\it smooth} hypersurfaces in projective space all MMM classes are trivial
in rational cohomology \cite{kn:RW}. That is, we need to go to a compactification
to see any interesting topology,  from this point of view. 
\item For sufficiently large $p$,  the direct images $\pi_{*}(K^{p})$ are
vector bundles over the moduli space, not just virtual bundles. As in other
moduli problems (see \cite{kn:AB} p. 582 for example), this gives restrictions
on their Chern characters and should lead to relations between the MMM classes.

\item Let $S$ be a smooth surface of general type and $\Pi$ be a 2-dimensional
subspace of $H^{0}(K_{S}^{q})$. If $q$ is sufficiently large and $\Pi$ is
generic this defines a Lefschetz pencil on $S$: the curves in the linear
system have at worst ordinary double points. This linear system gives a map
$$  \Gamma: \bC\bP^{1}= \bP(\Pi)\rightarrow \ocM_{g}, $$
where $2g-2= q(q+1)  c_{1}^{2}(S)$.
It seems possible that this construction could be developed to give some
connections between the putative enumerative theory of surfaces and curve-counting
invariants in the Deligne-Mumford spaces $\ocM_{g}$. (See the discussion
in 5.4  below.) 
\item To gain some insight into the obstruction  spaces for surface deformations
we consider smoothings of normal crossings. Let $S_{1}, S_{2}$ be surfaces
containing curves $C_{1}, C_{2}$ with normal bundles $N_{1}, N_{2}$. If $C_{1}$
and $C_{2}$ are isomorphic, say $C_{1}=C_{2}=C$, we form a singular surface
$\Sigma=S_{1}\cup_{C}S_{2}$ and study the smoothings  of this. The infinitesimal
deformations of the singularity correspond to sections of the line bundle
$N_{1}\otimes N_{2}$ over $C$.  We focus on the case when $S_{1}, S_{2}$
are cubic surfaces in $\bP^{3}$ and $C_{i}$ are in the linear system $\vert
\cO(3)\vert$ on $S_{i}$. In the end the smoothings we construct will be sextic
surfaces and from another point of view what we are studying is the degeneration
of smooth sextic surfaces to a union of two cubics. The general picture follows
the same pattern as the \lq\lq gluing'' techniques which have been employed
in all the Fredholm theories discussed in Section 1 (and there are important
recent developments in understanding the behaviour of K\"ahler-Einstein metrics
in such situations \cite{kn:Vetal}). 
The moduli space of cubic surfaces has dimension $4$ so the moduli space
of pairs $(S,C)$ of a cubic surface and curve in $\vert \cO(3)\vert$ has
dimension $22$. We need to study the matching problem for two pairs $(S_{1},C_{1}),
(S_{2}, C_{2})$ to have isomorphic curves $C_{1}, C_{2}$. The curves $C_{i}$
have genus $10$ and the moduli space $\cM_{10}$ has dimension $27$ so our
first guess is that the space of matching pairs has dimension $22+22-27=17$.
This will be true if the natural maps  between the various moduli spaces
have suitable transversality properties. The normal bundles $N_{i}$ have
degree $27$ and it follows from Riemann-Roch that, given a matching pair,
the dimension of $H^{0}(N_{1}\otimes N_{2})$ is $45$. So our first guess
is that the moduli space of the smoothed surface has dimension $45+17=62$
and this is indeed the {\it virtual} dimension of the moduli space of sextics,
as in item (1) above, but not the actual dimension which is $68$.

 The explanation for this deviation from the first-guess dimension count
is that the space of genus 10 curves which enter in the discussion is a lower-dimensional
subset $\cM_{10}'$ of $\cM_{10}$: the general curve of genus 10 cannot be
embedded as the intersection of two cubic surfaces in $\bP^{3}$. In fact
if $C=S_{1}\cap S_{2}$ where $S_{i}$ are cubic surfaces then $K_{C}=\cO(2)$
in other words $\cO(1)=K_{C}^{1/2}$ is a spin structure on $C$. At this point
we can refer to the article \cite{kn:Atiyah2} of Atiyah which
explains that the complex geometry of spin structures on complex curves can
be understood through the theory of skew-adjoint Fredholm operators. (An
operator $T$ on a complex Hilbert space is skew adjoint if $\langle x,Ty\rangle$
is a skew symmetric complex bilinear form.) The
$\db$-operator on $K_{C}^{1/2}$ becomes skew-adjoint with respect to standard
Hermitian structures. There are $2^{20}$ spin structures on a curve of genus
$10$ but for the purposes of this discussion, which only involves small deformations,
we can suppose that we have a distinguished one. Then the condition for a
curve $C$ to lie in $\cM'_{10}$ is that $H^{0}(K_{C}^{1/2})$ has dimension
$4$, while for generic curves it will have dimension $0$. Consider the abstract
situation of the space ${SFred}$ of skew-adjoint Fredholm operators on a
Hilbert space and the subset ${SFred}_{4}\subset { SFred}$ of operators with
kernel of dimension $4$. For $T\in { SFred}_{4}$ let $K_{T}$ be the 4-dimensional
kernel. Then one finds that normal bundle of ${SFred}_{4}$ in ${ SFred}$
has a $6$-dimensional fibre $\Lambda^{2} K_{T}$ at $T$. (This is a straightforward
extension of the obvious case when the Hilbert space has dimension $4$.)
In this way one finds that $\cM'_{10}$ has codimension $6$ in $\cM_{10}$
and the normal bundle is identified with $\Lambda^{2} H^{0}(K_{C}^{1/2})$.

The deformation theory of the singular surface $\Sigma$ yields a space ${\bf
T}^{1}(\Sigma)$ of infinitesimal deformations and an obstruction space ${\bf
T}^{2}(\Sigma)$. There is an exact sequence
$$ \dots H^{0}(N_{1}\otimes N_{2})\rightarrow {\bf T}^{1}(\Sigma)\rightarrow
{\bf T}^{1}(S_{1}, C_{1})\oplus {\bf T}^{1}(S_{2},C_{2})\rightarrow H^{1}(C;TC)\rightarrow
{\bf T}^{2}(\Sigma)\dots $$
where ${\bf T}^{1}(S_{i},C_{i})$ is the space of infinitesimal deformations
of the pair $(S_{i},C_{i})$. The term $H^{1}(C;TC)$ is the tangent space
of $\cM_{10}$ at $C$ but the incoming map in the sequence maps to the codimension
6 subspace $T\cM'_{10}$ and from the discussion above we get a map from $\Lambda^{2}\bC^{4}$
to ${\bf T}^{2}(\Sigma)$ which is in fact an isomorphism, fitting in with
what we saw in item (1) above. If we ran the whole discussion with $S_{1}$
a cubic surface but $S_{2}$ a quadric the corresponding curve $C$ has genus
$4$ and $K_{C}=\cO(1)$. The canonical system of a general curve of genus
$4$ embeds the curve in $\bP^{3}$ as the intersection of a cubic and a quadric
so our first-guess dimension count is correct in this case, fitting in with
the fact that $H^{2}(TS)$ vanishes for a smooth quintic surface. 

We see from this that the appearance of the obstruction spaces for surface
deformations is bound up with the special properties of curves on surfaces
stemming from the fact that curves are also divisors. This is well-known
in the enumerative  geometry of curves on surfaces. For a surface with $p_{g}>0$
({\it i.e.} $b^{+}>1$) \lq\lq most'' curves appear in families of the  wrong
dimension and do not contribute to the Gromov invariants. Our situation is
 different because a cubic surface has $p_{g}=0$. Start with our $4$-dimensional
family of complex structures on the smooth $4$-manifold underlying a cubic
surface and perturb this slightly to a $4$-dimensional family of almost-complex
structures. Then we still have a $22$-dimensional family of pairs $(S,C)$
but now we expect that the matching problem will behave in a generic way
and that the space of matching pairs $(S_{1}, C_{1}), (S_{2}, C_{2})$ will
have   dimension 17, rather than 23 as occurs in the integrable case. (Strictly
we should pass to real dimensions here, since the moduli spaces will not
have complex structures.)

\end{enumerate}

\section{Study of an example}

For most surfaces that one can construct easily the moduli spaces have very
large dimension. Imposing symmetry by a finite group allows us to cut down
the dimension to get manageable spaces but still exhibiting some essential
features of the situation. In this Section we will discuss one example of
this kind and compute the  enumerative invariants (modulo some
foundational assumptions). 

\

{\bf Remark} The significance of the numerical invariants $a=c_{1}^{2},b=
\chi$ of a surface
is that they define  the Hilbert polynomial for ${\rm dim}\ H^{0}(K_{S}^{p})$.
Fix a finite group $\Gamma$ and consider surfaces of general type with a
$\Gamma$ action. The spaces $H^{0}(K_{S}^{p})$ are then representations of
$\Gamma$ so the Hilbert polynomial extends to a \lq\lq $\Gamma$-Hilbert function'',
taking values in the representation ring of $\Gamma$ which is the
generalisation of the pair $(a,b)$.   

\subsection{A family of sextic surfaces}
Take coordinates $(x_{1}, y_{1}, x_{2}, y_{2})$ on $\bC^{4}$ and let $\zeta$
be a primitive sixth root of unity. Let $G$ be the subgroup of $GL(4,\bC)$
generated by:
$$ (x_{1}, y_{1}, x_{2}, y_{2})\mapsto (\zeta x_{1},\zeta^{-1} y_{1}, x_{2},
y_{2})$$
$$(x_{1}, y_{1}, x_{2}, y_{2})\mapsto (x_{1}, y_{1}, \zeta x_{2},\zeta^{-1}
y_{2}) $$
$$ (x_{1}, y_{1}, x_{2}, y_{2}) \mapsto (x_{2}, y_{2}, x_{1}, y_{1}). $$

Let $V$ be the vector space of polynomials of degree $6$ invariant under
$G$. These have the form
\begin{equation}  \alpha x_{1}^{6}+ \beta y_{1}^{6}+\alpha x_{2}^{6}+ \beta
y_{2}^{6} + A Q_{+}^{3} + B Q_{+} Q_{-}^{2} , \end{equation}
where $Q_{\pm}= x_{1} y_{1} \pm x_{2} y_{2}$. The $\bC^{*}$-action on $\bC^{4}$
$$  (x_{1}, y_{1}, x_{2}, y_{2})\mapsto (\lambda x_{1}, \lambda^{-1} y_{1},\lambda
x_{2},\lambda^{-1} y_{2}), $$
induces an action on $V$:
$$  (\alpha,\beta_, A, B)\mapsto (\lambda^{6}\alpha, \lambda^{-6}\beta, A,B).
$$
The stable points in $V$ for the torus action  are those where $\alpha$ and
$ \beta$ are non-zero so each stable orbit in $V$ contains a  representative
$$  \alpha (x_{1}^{6}+ y_{1}^{6}+ x_{2}^{6}+ y_{2}^{6}+) + A Q_{+}^{3}+ B
Q_{+} Q_{-}^{2}, $$
which is unique up to change in the sign of $\alpha$.  Thus we get a moduli
space $M$ of \lq\lq GIT stable'' sextic surfaces with this $G$-action which
is the quotient of $\bC^{2}$ by $\pm 1$,  where a point $(A,B)$ in $\bC^{2}$
corresponds to the surface
$$  S_{AB}= \{ x_{1}^{6}+y_{1}^{6}+ x_{2}^{6}+ y_{2}^{6}+ A Q_{+}^{3}+ B
Q_{+} Q_{-}^{2}=0\}\subset \bC\bP^{3}. $$

The locus $\Delta\subset \bC^{2}$ of points $(A,B)$ where $S_{AB}$ is singular
has four components. Write $F=(3A-B)$ and for  $\epsilon_{1},\epsilon_{2}\in
\{\pm 1\}$  write $G_{1}=3A+3B+5 \epsilon_{1}, G_{2}=3A+3B+5 \epsilon_{2}$.
Then for each of these $4$ choices of signs there is a component of $\Delta$
with equation
$$ (F^{2}- G_{1}G_{2})^{2}= 4 F^{2}(G_{1}-F)(G_{2}-F). $$
The singularities that arise are mild and will not enter  further in our
discussion.

The reason for choosing this  example is that the virtual dimension of $M$
is not equal to the actual dimension. In general if we have a surface $S$
with the action of a finite group $\Gamma$ the deformation theory, for surfaces
with $\Gamma$-action, works in the obvious way. The group $\Gamma$ acts on
$H^{i}(TS)$, infinitesimal deformations are given by the invariant part $H^{1}(TS)^{\Gamma}$
and obstructions in $H^{2}(TS)^{\Gamma}$. We go back to amplify  the  description
of the obstruction spaces for hypersurfaces we mentioned in Section 3, for
the case of sextic surfaces. On $\bP^{3}= \bP(U)$ we have the dual Euler
sequence 
$$  T^{*}\bP^{3}(1)\rightarrow \underline{U^{*}}\rightarrow \cO(1), $$ and
taking
tensor product with $\cO(1)$ gives

$$ T^{*}\bP^{3}(2)\rightarrow \underline{U^{*}}(1)\rightarrow \cO(2). $$
Taking sections we get
$$   H^{0}(T^{*}\bP^{3}(2))\rightarrow U^{*}\otimes U^{*}\rightarrow s^{2}(U^{*}),
$$
which shows that $H^{0}(T^{*}\bP^{3}(2))$ is canonically isomorphic to $\Lambda^{2}U^{*}$.
Now for a smooth sextic surface $S\subset \bP^{3}$ cut out by a section
$s$ of $\cO(6)$ we have a restriction map $T^{*}\bP^{3}\vert_{S}\rightarrow
T^{*}S$ and an isomorphism of line bundles $K_{S}=\cO(2)$. The second isomorphism
depends, up to a factor, on the choice of $s$ and a volume form $\Lambda^{4}U=\bC$.
Combining these ingredients, we get a map $r: \Lambda^{2}U^{*}\rightarrow
H^{0}(T^{*}S\otimes K_{S})$ which one readily sees is an isomorphism. Now
suppose that a group $\Gamma$ acts on $V$, preserving $S$. If $\Gamma$ acts
with determinant $1$ on $V$ and if $\Gamma$ preserves the section $s$ cutting
out $S$ then it follows that $r$ is a $\Gamma$-equivariant map, for the standard
actions on source and target. This holds in our situation, with $\Gamma=G$,
and we conclude  that the dual of the obstruction space for our surfaces
with $G$-action is the $G$-invariant part of $\Lambda^{2}U^{*}$. One finds
that this is $1$-dimensional, spanned by
$  \omega= dx_{1} dy_{1} + dx_{2} dy_{2} $.
So the virtual (complex) dimension of $M$ is $1$. 

\

There is an obvious \lq\lq naive'' compactification $\oM_{GIT}$ of $M$, which
is to take $\bC^{2}\subset \bC\bP^{2}$ and the quotient of $\bC\bP^{2}$ by
$\pm 1$. This is the Geometric Invariant theory compactification obtained
by adding \lq\lq polystable'' points. The points at infinity in $\oM_{GIT}$
correspond to solutions of the equation $A Q_{+}^{3}+ BQ_{+} Q_{-}^{2}=0$.
If $A,B\neq0$ these form a union of three quadrics meeting in four lines;
if $B=0$ we get a quadric $\{Q_{+}=0\}$ with multiplicity $3$  and if $A=0$
a union of one  quadric $\{Q_{+}=0\}$ and another $\{Q_{-}=0\}$ with multiplicity
$2$. None of these objects is allowed in the KSBA compactification. 

 \

The naive compactification $\oM_{GIT}$ is a toric surface and has a toric
description $(P,\Lambda)$ where $\Lambda\subset \bZ^{2}$ is the lattice of
pairs $(n,m)\in \bZ^{2}$ with $n+m$ even and $P$ is the triangle with vertices
$(0,0), (-6,0), (0,-6)$. (Of course, the complex structure only determines
the \lq\lq fan'' of normals to the edges of $P$ and we could scale $P$ by
a factor.) This toric structure is partly accidental---the torus does not
act on the family of surfaces parametrised by the moduli space. 

\

\

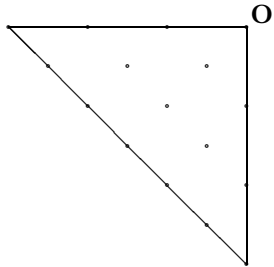
\setlength{\unitlength}{.3pt}

\begin{picture}(320,320)
\put(310,10){\line(0,1){300}}
\put(10,310){\line(1,0){300}}
\put(310,10){\line(-1,1){300}}
\put(210,310){\circle{4}}
\put(110,310){\circle{4}}
\put(10,310){\circle{4}}
\put(310,210){\circle{4}}
\put(310,110){\circle{4}}
\put(310,10){\circle{4}}
\put(260,260){\circle{4}}
\put(160,260){\circle{4}}
\put(260,160){\circle{4}}
\put(210,210){\circle{4}}
\put(260,60){\circle{4}}
\put(210,110){\circle{4}}
\put(110,210){\circle{4}}
\put(60,260){\circle{4}}
\put(160,160){\circle{4}}
\put(310,310){\circle{4}}
\put(315,315){{\bf O}}
\end{picture}

\subsection{The KSBA compactification}

In this subsection we will describe another compactification $\oM$ of the
moduli space $M$. Before going to the general picture we begin by discussing
the $1$-parameter family of surfaces $S_{A,0}$ which is more straightforward.
The basic  point is that the correct limit when $A\rightarrow \infty$ is
a triple branched cover of the quadric $Q_{+}=0$. Let $\pi:Y\rightarrow \bC\bP^{3}$
be the  cyclic triple cover of $\bC\bP^{3}$ branched over the smooth surface
$S_{0,0}$.  If $s$ is the section of $\cO(6)$ over $\bC\bP^{3}$ cutting out
$S_{0,0}$ there is by construction a section $\eta$ of $\pi^{*}(\cO(2))$
over $Z$ with $\eta^{3}= s$. We  also have another section $Q_{+}$ of $\pi^{*}(\cO(2))$.
Let $W$ be the hypersurface in $Z\times \bP^{1}$ defined by the equation
written
$\eta=\lambda Q_{+}$ in terms of an affine coordinate $\lambda$ on $\bP^{1}$.
So we have a projection $W\rightarrow \bP^{1}$ and the fibre over a finite
point $\lambda$ can be identified with $S_{A,0}$ where $A=\lambda^{3}$. The
fibre over infinity is a smooth surface $S_{III}$: the triple cover of the
quadric $Q_{+}=0$ branched over the intersection with $S_{0,0}$. The fact
that we have to take the cube root $\lambda=A^{1/3}$ to construct the family
is the usual orbifold phenomenon: the $\bZ/3$-action of the triple cover
means that $S_{III}$ has a larger automorphism group than the generic surface
$S_{A,0}$.

\

Our compactification $\oM$ is also a toric variety. Let $\Pi$ be the quadrilateral
with vertices $0=(0,0), II=(0,-6), III=(-4,0), IV=(-2,-6)$ and $\Lambda\subset
\bZ^{2}$ be the same lattice as above. Then we define $\oM$ to be  polarised
  toric variety corresponding to $(\Pi,\Lambda)$ (we discuss the polarisation
later).  To see this as a moduli
space we want to take a toric chart for each  vertex,  construct a corresponding
family of surfaces and glue these together over the intersections of the
charts. For brevity we will  do this in full for just two charts corresponding
to the vertices $0, IV$. The first we already have: it is the family of surfaces
$S_{A,B}$ parametrised by $\bC^{2}/\pm 1$ and the centre of the chart corresponds
to the surface $S_{O}=S_{0,0}$.

\

\

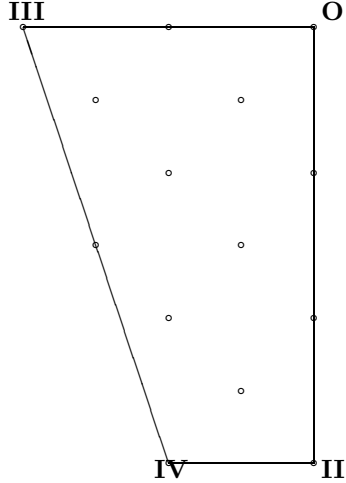
\setlength{\unitlength}{.55pt}
\begin{picture}(320,320)
\put(110,10){\line(1,0){100}}
\put(110,10){\line(-1,3){100}}
\put(10,310){\line(1,0){200}}
\put(210,10){\line(0,1){300}}
\put(100,0){{\bf IV}}
\put(215,0){{\bf II}}
\put(0,315){{\bf III}}
\put(210,210){\circle{4}}
\put(110,210){\circle{4}}
\put(210,110){\circle{4}}
\put(160,260){\circle{4}}
\put(110,110){\circle{4}}
\put(210,10){\circle{4}}
\put(110,10){\circle{4}}
\put(110,310){\circle{4}}
\put(160,160){\circle{4}}
\put(160,60){\circle{4}}
\put(60,260){\circle{4}}
\put(60,160){\circle{4}}
\put(210,310){\circle{4}}
\put(10,310){\circle{4}}
\put(215,315){{\bf O}}
\end{picture}

\

\

Before going on to the vertex $IV$ we describe the surfaces corresponding
to the rest of the boundary of the quadrilateral. Recall that a smooth quadric
 in $\bC\bP^{3}$ can be identified with $\bP^{1}\times\bP^{1}$ and take standard
affine co-ordinates $(s,t)$ on the latter,  so we have four lines $s=0,s=\infty,
t=0,t=\infty$. This configuration of four lines will appear often in what
follows and we will denote it by $\Lambda$.  Let
$C_{\mu}\subset \Sigma$ be the curve in the linear system $\cO(6,6)$ over
$\bP^{1}\times \bP^{1}$ with affine equation
\begin{equation} 1 + s^{6}+t^{6}+ s^{6}t^{6}+ \mu s^{3} t^{3} =0\end{equation}
\begin{itemize}
\item The segment from $O$ to $III$ corresponds to the surfaces $S_{A0}$,
and the segment from $O$ to $II$ to the $S_{0B}$.
\item The open segment from $III$ to $IV$ corresponds to triple covers of
$\bP^{1}\times \bP^{1}$ with simple branching over the curves $C_{\mu}, C_{-\mu}$,
for $\mu\in
\bC\setminus \{0\}$. The limit as $\mu\rightarrow 0$ is the surface $S_{III}$
we discussed above: a $\bZ/3$-cyclic cover branched over $C_{0}$. 
\item The  segment from $II$ to $IV$, including the end-point $II$ but not
$IV$, corresponds to the following family of surfaces. For $\mu\in \bC$ take
the double cover of $\bP^{1}\times \bP^{1}$ branched over the divisor in
$\cO(8,8)$ given by the union of $C_{\mu}$ and 
the four lines $\Lambda$. So $\Lambda$ lifts to the double cover.  Each line
meets $C_{\mu}$ in $6$ points and these give ordinary double points in the
branched cover. Blow up these $24$ points and then collapse the proper transform
of $\Lambda$ to a point. Taking $\mu=0$ gives the surface $S_{II}$ corresponding
to the vertex $II$. 
\end{itemize}

\

To study the vertex $IV$ we let $\bP^{5}_{w}$ be the weighted projective
space with homogeneous coordinates $(x_{1}, y_{1}, x_{2}, y_{2}, h_{+}, h_{-})$
and weights $1$ on the first four coordinates and $2$ on the last two. We
consider $\bP^{3}$ as embedded in $\bP^{5}_{w}$ in the obvious way and write
$\bP^{1}_{h}$ for the line $\{(0,0,0,0, h_{+}, h_{-})\}$ in $\bP^{5}_{w}$.
For parameters $\alpha,\beta\in \bC$, let $S^{\alpha,\beta}$ be the complete
intersection in $\bP^{5}_{w}$ defined by the equations:
$$   x_{1}^{6}+y_{1}^{6}+ x_{2}^{6}+y_{2}^{6}+ h_{+}^{3}+ h_{+} h_{-}^{2}
= 0, $$
$$ x_{1}y_{1}= \alpha h_{+} + \beta h_{-}, $$
$$ x_{2}y_{2}= \alpha h_{+}- \beta h_{-}. $$
If $\alpha, \beta$ are both nonzero we can write $h_{+}=(x_{1}y_{1}+x_{2}y_{2})/2\alpha
$ and $h_{-}= (x_{1}y_{1}-x_{2}y_{2})/2\beta$  and, substituting into the
first equation.  we get
$$ 8 ( x_{1}^{6}+y_{1}^{6}+ x_{2}^{6}+y_{2}^{6}) + \alpha^{-3} (x_{1} y_{1}+
x_{2}y_{2})^{3}+ \alpha^{-1}\beta^{-2}(x_{1}y_{1}+ x_{2} y_{2})(x_{1}y_{1}-
x_{2}y_{2})^{2}=0. $$

The surface $S^{\alpha,\beta}$ does not meet the line $\bP^{1}_{h}$ in $\bP^{5}_{w}$
and projection from this line maps $S^{\alpha,\beta}$ to $S_{A,B}$ in $\bC\bP^{3}$
where
\begin{equation} 8 A= \alpha^{-3} \ \ \ \ 8 B= \alpha^{-1}\beta^{-2} \end{equation}

We define $S_{IV}$ to be the surface $S^{0,0}$. Let $\bP^{4}_{w}$ be the
weighted projective space with co-ordinates $(x_{i}, y_{j}, h_{-})$. The
projection
from $S_{IV}$ to $\bP^{4}_{w}$ is well defined and the image is a cone $\Cone(\Lambda)$
over the configuration $\Lambda$ of $4$-lines in $\bP^{3}$. The projection
exhibits $S_{IV}$ as a triple branched cover of $\Cone(\Lambda)$. The surface
$S^{0,0}$ meets $\bP^{1}_{h}$ in three points $h_{+}=0, h_{+}=h_{-}, h_{+}=-h_{-}$
and the triple cover maps these  three points  to the vertex of the cone.
 
  Next consider the case when $\alpha$ is zero but $\beta$ is not. Then the
equations give
$x_{1}y_{1}+x_{2}y_{2}=0$ and the projection from $S^{0,\beta}$ to $\bP^{3}$
has image this quadric surface. Take affine co-ordinates $x_{1}=s, y_{1}=t,
x_{2}=st, y_{2}=1$ on this quadric, as before. We have $h_{-}= \beta^{-1}
x_{1} y_{1}= \beta^{-1} st$ and our surface is defined by the equation
$$   (1+ s^{6}+t^{6}+ s^{6}t^{6}) + h_{+}^{3}+ \beta^{-2}s^{2} t^{2}h_{+}
=0. $$ The projection to the quadric exhibits $S^{0,\beta}$ as a triple cover
with branch locus 
$$  (1+ s^{6}+t^{6}+ s^{6}t^{6})\pm \frac{2}{3\beta^{3}\sqrt{-3}} s^{3}t^{3}=0,
$$
which agrees with our previous discussion of the boundary segment from $III$
to $IV$.

  The case  $\beta=0$ is similar but a little more complicated. The equations
give $x_{1}y_{1}-x_{2}y_{2}=0$ defining another smooth quadric in $\bC\bP^{3}$
containing the same line configuration $\Lambda$. We  parametrise this quadric
surface by $x_{1}=s, y_{1}=t, x_{2}= st, y_{2}=-1$. The surface $S^{\alpha,0}$
meets the line $\bP^{1}_{h}$ in the point $h_{+}=0$ and the projection from
$S^{\alpha,0}$ to $\bP^{3}$ is not defined there. Blowing up this point we
get a surface  $\tilde{S}^{\alpha,0}$ which  maps to $\bP^{3}$ with image
the quadric. Writing $h_{+}=\alpha^{-1} x_{1}y_{1}$ this blown up surface
$\tilde{S}^{\alpha,0}$ is defined by the equation
$$   (1+ s^{6}+t^{6}+ s^{6}t^{6}) + \alpha^{-3} s^{3}t^{3} + \alpha^{-1}st
h^{2}_{-}=0$$

If $s=0$ and $t^{6}=-1$ the co-ordinate $h_{-}$ is unconstrained and we get
a line in $\tilde{S}^{\alpha,0}$; similarly for $s=\infty, t=0,\infty$. Collapsing
these 24 lines in $\tilde{S}^{\alpha,0}$ gives the double cover of the quadric
branched over the $\cO(8,8)$ divisor
$$  st  \left( (1+ s^{6}+t^{6}+ s^{6}t^{6}) + \alpha^{-3} s^{3}t^{3}\right)=0.
$$
which agrees with the previous discussion for the boundary
segment from $II$ to $IV$.

From another point of view, there is a well-known toric degeneration of the
quadric $\bP^{1}\times \bP^{1}$ to the cone $\Cone(\Lambda)$. We start with
$S_{IV}$, the triple cover of the cone, and deform this to  triple covers
of
$\bP^{1}\times\bP^{1}$ to get the surfaces $S^{0,\beta}$, which are smooth
for small $\beta$.  Similarly for the $S^{\alpha,0}$ but with the extra complication
due to the point $h_{+}=0$ in $\bP^{1}_{h}$ which is not smoothed by the
deformation and remains a singular point in $S^{\alpha,0}$.

We can make  similar constructions in toric charts corresponding to the vertices
$II,III$ but
the key formula for constructing the toric moduli space $\oM$ is (11). First,
we see that only $\beta^{2}$ appears, so we write $\beta^{2}=\gamma$ and
we have $A^{-1} = \alpha^{3}, B^{-1}= \alpha\gamma$. Recall that $(-A,-B)$
defines the same point in the moduli space as $(A,B)$. Thus changing the
sign of $\alpha$ does not change the surface $S^{\alpha,\beta}$. Consider
a monomial $A^{-p}B^{-q}$ on $(\bC^{*})^{2}$ which is equal to $\alpha^{p+3q}\gamma^{q}$.
For this to descend to a well-defined function on the moduli space (i.e.
invariant under change of  sign of $\alpha$) we need $p+3q$ to be even, so
$p+q$ is even. For the monomial to extend holomorphically over $\alpha=\gamma=0$
we need $3p+q\geq 0, q\geq 0$. So the holomorphic functions on a neighbourhood
of the point corresponding to $S_{IV}$ in the moduli space have a basis given
by the intersection of our lattice  $\Lambda\subset\bZ^{2}$ with the convex
set $\{(p,q)\in \bR^{2}: 3p+q\geq 0, q\geq 0\}$. Standard toric theory shows
that  the compact space $\oM$ is defined by the polytope given by intersecting
further with a set $\{(p,q): p\leq C_{1}, q\leq C_{2}\}$ for any fixed $C_{1},
C_{2}>0$. We have taken $C_{1}=2, C_{2}= 6$ and then translated the quadrilateral
so that the origin is at the vertex $O$.

The author has a strong feeling that this moduli space $\oM$ is the KSBA
moduli space (for these surfaces with $G$-action) but he is not qualified
to certify that as a definite fact. In any case we will proceed with out
study based on that assumption. 

\

{\bf Remarks}

 1. We can make the same constructions in $\bP^{9}$ using the canonical
embeddings of our surfaces $S_{A,B}$, but then we are in high  codimension
with many equations and variables  which do not play any real role. The advantage
of the weighted projective space is that it allows us to bring in just the
two sections of $K_{S}$ which are really relevant. It seems likely to the
author that the moduli space we are constructing is the Chow stable moduli
space, under the canonical embedding. 

2. Starting with the GIT compactification $\oM_{GIT}$ we can get the compactification
$\oM$ by performing weighted blow-ups at the two points corresponding to
the vertices $ (-6,0),(0,-6)$ and contracting the proper transform of the
line at infinity. For GIT moduli spaces there are  techniques of Jeffrey
and Kirwan \cite{kn:JK} which, in favourable cases, can be used to calculate
pairings of the kind we are concerned with, so the comparison of the different
compactifications is a relevant topic. Laza investigates this comparison,
for another moduli problem, in \cite{kn:Laza}.

\subsection{Calculations in cohomology}

The moduli space $\oM$ has virtual complex dimension $1$ so there is a virtual
fundamental class $\zeta \in H_{2}(\oM)$. We have two MMM classes, associated
to the characteristic classes $c_{1}^{3}$ and $c_{1}c_{2}$. The goal of this
subsection is to calculate the pairing of $\zeta$ with these two classes.
By standard toric theory $H^{2}(\oM;\bQ)$ is two dimensional. Each edge of
the quadrilateral $\Pi$ corresponds to a $2$-sphere in $\oM$ and so defines
a homology class. Let $D_{II}, D_{III}$ be the 2-spheres corresponding to
the edges from $O$ to $III$ and $O$ to $II$ respectively and use the same
symbols for their homology classes. These give a basis for $H_{2}(\oM,\bQ)$
in which we will do our calculations. We begin by calculating the $2$-dimensional
MMM classes.  

Consider first a general situation. Let $Y$ be a Calabi-Yau 3-fold and $L\rightarrow
Y$ an ample line bundle. Serre duality implies that the Hilbert polynomial
of $Y$ is odd, say:
\begin{equation}   {\rm dim}\  H^{0}(L^{p})= H(p)= \alpha p^{3} + \beta p,
\end{equation}
for sufficiently large $p$ (in fact $p\geq 1$). Suppose that  $s_{0}, s_{1}$
are sections of $L$ defining a Lefschetz pencil on $Y$. Then we have a subvariety
$W\subset Y\times \bP^{1}$ cut out by the section of $L\otimes \cO(1)$ written
in an affine coordinate on $\bP^{1}$ as $s_{0}- \lambda s_{1}$. Write $\pi:W\rightarrow
\bP^{1}$ for the projection and
$K_{v}$ for the relative canonical bundle. The adjunction formula gives $K_{v}=
L\otimes \cO(1)$. For large enough $p$ we have a vector bundle
$ V_{p}= \pi_{*}(K_{v}^{p}) $
over $\bP^{1}$ of rank $r_{p}$ and degree $d_{p}$. We want to find $d_{p}$
in terms of the  data $\alpha,\beta$. We could apply the Riemann-Roch theorem
for families but in this situation there is a more elementary direct argument.
If $V$ is a vector bundle over $\bP^{1}$ of rank $r$ and degree $d$ then
for large enough $q$ we have
$$  {\rm dim}\  H^{0}(\bP^{1}, V\otimes \cO(q))= q r + (d+r). $$
Applied to $V_{p}$ we get
$$ {\rm dim}\  H^{0}(W, L^{p}\otimes \cO(p+q))= q r_{p} + (d_{p}+r_{p}),
$$
for large enough $p,q$. The  restriction sequence for $W\subset M\times \bP^{1}$
gives
$$  {\rm dim}\  H^{0}(W,L^{p}\otimes \cO(p+q))= H(p)(p+q+1)- H(p-1)(p+q),
$$
and, comparing the two formulae, we get
$$ d_{p}= p H(p)- (p-1) H(p-1). $$
Writing $\up=p-1/2$ as before, this is
\begin{equation}  d_{p}= \alpha \left( 4\up^{3}+ \up \right) + 2\beta \up.
\end{equation}
 Comparing with the formulae (6),(7) for the Knudsen-Mumford line bundles
$\cL_{0},\cL_{2}$ we see that
\begin{equation} \langle c_{1}(\cL_{0}), [\bP^{1}]\rangle= -(\alpha+\beta)
\ ,\ \ \langle c_{1}(\cL_{2}), [\bP^{1}]\rangle= 12 \alpha. \end{equation}

To apply this in our situation we begin with the 2-sphere $D_{III}\subset
\oM$. The family that we discussed at the beginning of subsection 6.2 above
can be embedded in weighted projective space $\bP^{4}_{w}$. In our co-ordinates
$(x_{1}, y_{1}, x_{2}, y_{2}, h_{+})$ we let $Y$ be the degree $6$ hypersurface
defined by the equation
$$   x_{1}^{6}+y_{1}^{6}+x_{2}^{6}+ y_{2}^{6}+  h_{+}^{3}=0. $$    
The adjunction formula in $\bP^{4}_{w}$ shows that $Y$ is a Calabi-Yau 3-fold.
We take the line bundle $\cO(2)$ over $Y$ and consider the pencil $h_{+}-
\lambda (x_{1}y_{1}+ x_{2}y_{2})$, so for finite $\lambda$ the fibre of
$\pi:W\rightarrow\bP^{1}$ is the surface $S_{A,0}$ with $A=\lambda^{3}$ and
the fibre over $\infty$ is $S_{III}$. To find the Hilbert polynomial $H(p)$
for this pair $(Y,L)$ let $n_{p}$ be the dimension of the space of homogeenous
polynomials of degree $2p$ on $\bC^{4}$, so
$$  n_{p}= \frac{1}{6} (2p+1)(2p+2)(2p+3). $$
Then if we write $N_{p}= {\rm dim} H^{0}(\bP^{4}_{w}, \cO(2p))$ we have
$$  N_{p}= n_{p}+ n_{p-1}\dots + n_{0}, $$
whereas the restriction sequence for $Y\subset \bP^{4}_{w}$ gives, for large
enough $p$,
$   H(p)= N(p)- N(p-3)$.  So we conclude that
$ H(p)= n_{p}+ n_{p-1}+n_{p-2}$ and this gives
\begin{equation}  H(p)= 4p^{3}+ 7p. \end{equation}

Since the true parameter on $D_{III}$ is $A^{2}= \lambda^{6}$ we get a factor
of $1/6$ in the formulae, and we arrive at
\begin{equation} \langle c_{1}(\cL_{0}), [D_{III}]\rangle= -11/6 \ \ ,\ \
\langle c_{1}(\cL_{2}), [D_{III}]\rangle= 8. \end{equation}

The discussion for $D_{II}$ is very similar. This time we define a degree
$6$ hypersurface  $Y'\subset \bP^{4}_{w}$ by the equation
$$     x_{1}^{6}+y_{1}^{6}+x_{2}^{6}+ y_{2}^{6}+ (x_{1} y_{1}+ x_{2} y_{2})
h_{-}^{2}=0, $$  
and consider the pencil $(x_{1}y_{1}- x_{2} y_{2}) - \lambda h_{-}$. The
only
difference is that $Y'$ is singular, with a singular point at $P_{\infty}=(0,0,0,0,1)$
but the singularity does not affect the calculations. The true parameter
on $D_{II}$ is $B^{2}= \lambda^{4}$ so we get
\begin{equation} \langle c_{1}(\cL_{0}), [D_{II}]\rangle= -11/4 \ \ ,\ \
\langle c_{1}(\cL_{2}), [D_{II}]]\rangle= 12. \end{equation}
 and we see that $c_{1}(\cL_{0})= -(11/48) c_{1}(\cL_{2})$ in $H^{2}(\oM,
\bQ)$.

The quadrilateral $\Pi$ has been chosen to correspond to the polarisation
$\Omega= c_{1}(\cL_{2})/4\in H^{2}(\oM)$. By general toric theory
$$  \langle \Omega^{2}, [\oM]\rangle = 2 {\rm Area}_{\Lambda}\ (\Pi), $$
where ${\rm Area}_{\Lambda}$ means the area relative to the lattice, which
is half the standard area. So we see that
$  \langle \Omega^{2}, [\oM]\rangle = 18$ and hence
 
\begin{equation}   \langle c_{1}(\cL_{2})^{2}, [\oM] \rangle= 16.18=288.
\end{equation}

We now consider the virtual fundamental class $\zeta\in H_{2}(\oM;\bQ)$.
Due to the foundational gap expressed in our \lq\lq Premise'' we do not have
an  official definition of this,  but we calculate in what seems the appropriate
way. Recall that the standard symplectic form $\omega$ on $\bC^{4}$ defines
a section $s_{\omega}$ of $T^{*}\bP^{3}(2)$.  This section has no zeros so
defines a rank-2 sub-bundle
$E\subset T^{*}\bP^{3}$. (This is a well-known object: it is the \lq\lq null
correlation bundle'' which corresponds via twistor theory to the standard
Yang-Mills instanton on $S^{4}$.) The sub-bundle $E$ is a holomorphic contact
structure on $\bC\bP^{3}$ so there is no surface $\Sigma\subset \bC\bP^{3}$
(even locally) such that the restriction of $s_{\omega}$ to $T^{*}\Sigma(2)$
vanishes. We lift $s_{\omega}$ by the projection $\bP^{4}_{w}\setminus P_{\infty}
\rightarrow \bP^{3}$ to define a section $s'_{\omega}$ of $T^{*}\bP^{4}_{w}(2)$
away from the point $P_{\infty}$. Now consider the degree $6$ hypersurface
$Y\subset \bP^{4}_{w}$ and pencil in $\cO(2)$ defining a family $\pi:W\rightarrow
\bP^{1}$ as above,  with $W\subset Y'\times \bP^{1}\subset \bP^{4}_{w}\times
\bP^{1}$. Lifting the section $s'_{\omega}$ to the product and restricting
to $W$ we finally get a section $s''_{\omega}$ of $T^{*}W \otimes K_{v}\otimes
\cO(-1)$. Let $\Delta\subset \bC\subset \bP^{1}$ be the finite set corresponding
to singular surfaces in the pencil. We have a line bundle $\cL\rightarrow
\bC\setminus \Delta$ with fibres the $G$-invariant part of the obstruction
spaces $H^{2}(TS)$ and a dual line bundle $\cL^{*}$ with fibres the $G$-invariant
part of $H^{0}(T^{*}S\otimes K_{S})$. The section $s''_{\omega}$ defines
by restriction to fibres a non-vanishing section of $\cL^{*}(-1)$. On the
smooth part of a singular fibre or on the fibre over $\lambda=\infty$ the
section restricts to a finite and non-vanishing section of $T^{*}S\otimes
K_{S}$, so the natural extension of $\cL^{*}$ to $\bP^{1}$ is isomorphic
to $\cO(1)$. The same discussion applies to the family corresponding to $D_{II}$:
we just restrict to the smooth part.  Taking account again of the coverings
we get
\begin{equation}   \zeta. D_{III} =-1/6  \  \ \ \ \  \zeta. D_{II}=- 1/4
. \end{equation}
 
 \

 We make a short digression to consider further the $\bZ/3$-cyclic cover
$S_{III}$ which gives insight into the denominators in these formulae. Let
$V_{+}, V_{-}$
be $2$-dimensional complex vector spaces with fixed isomorphisms $\Lambda^{2}V_{\pm}=\bC$
and write $\bP^{1}\times \bP^{1}= \bP(V_{+})\times \bP(V_{-})$. Let $S\rightarrow
\bP^{1}\times \bP^{1}$ be a $\bZ/3$-cyclic cover branched over a smooth curve
of bi-degree $(6,6)$. Using standard theory one finds that
\begin{equation}   H^{0}(T^{*}S\otimes K_{S}) = s^{2}(V^{*}_{+}) \oplus s^{2}(V^{*}_{-})\end{equation}
 Now let $U$ be the $4$-dimensional space $V_{1}\otimes V_{2}$. The trivialisations
of $\Lambda^{2}V_{\pm}$ define a nondegenerate quadratic form on $U$ and
$\bP^{1}\times \bP^{1}$ is embedded in $\bP(U)$ as the corresponding quadric
surface. This is the usual description of $4$-dimensional (complex) oriented
Euclidean geometry in terms of spin spaces $V_{\pm}$.We have the usual splitting
of the $2$-forms
$  \Lambda^{2} U^{*}= \Lambda^{2}_{+}\oplus \Lambda^{2}_{-}$ and one has
the spinor description
$  \Lambda^{2}_{\pm}= s^{2}(V^{*}_{\pm})$.
 What we essentially see from this is that the obstruction spaces behave
in a simple way for a family of smooth sextic surfaces converging to a triple
cover of a quadric. In the case at hand the group $G$ acts on $V_{\pm}$ and
$G$-invariant piece corresponding to the self-dual form $\omega$ matches
up with a $1$-dimensional invariant subspace in $s^{2}(V^{*}_{+})$. However
the
isomorphism (20) depends on an isomorphism $K_{S}=\cO(2,2)$ and the $\bZ/3$
covering group acts non-trivially on this. So on our $\bP^{1}$ covering $D_{III}$,
with parameter $\lambda$, we have a line bundle $\cL=\cO(-1)$ but the additional
$\bZ/3$ symmetry group of $S_{III}$ acts non-trivially on the fibre over
$\infty$. 

Similarly there is an additional $\bZ/2$ symmetry group of the surface $S_{O}$
which works the equivalence $S_{A,B}=S_{-A,-B}$. This is given by
$x_{i}\mapsto -x_{i}, y_{j}\mapsto y_{j}$ and so takes $\omega$ to $-\omega$.
So for both $S_{O}$ and $S_{III}$ if we take account of their full symmetry
groups the corresponding invariant part of the obstruction space vanishes.
This means that under a generic perturbation both $S_{O}$ and $S_{III}$ \lq\lq
persist'': more precisely, slightly perturbed versions of them, with the
additional symmetries,  persist.  The picture is that  the solutions $Z$
of
the perturbed problem can be viewed (approximately) as a subset of $\oM$
but the symmetries of the situation force this subset to contain the points
in $\oM$ corresponding to  vertices of the quadrilateral. The moduli
space  $\oM$ is only a rational homology manifold at these points and the
intersection pairing on integral homology takes rational values. 

\

To sum up: we have three classes in $H^{2}(\oM, \bQ)$ given by $c_{1}(\cL_{0}),
c_{1}(\cL_{2})$ and the Poincar\'e dual of $\zeta$ and these are are all
proportional:
$$  PD[\zeta]= (1/48) c_{1}(\cL_{2})= -(1/11)c_{1}(\cL_{0}).$$  
Using (18) we find the  pairings
\begin{equation}  \langle \zeta, c_{1}(\cL_{2})\rangle = -6 \ \ \ ,\ \ \
\langle \zeta, c_{1}(\cL_{0})\rangle=  11/8 . \end{equation}
Using the formulae (8) we see that the pairings with $I(c_{1}^{3}), I(c_{1}c_{2})$
 are $12$ and $-11/192$ respectively. 
 
\subsection{Curves}

Let $\eta$ be a 12th. root of unity and take the  action of $\bZ/12$ on $\bC^{2}$
generated by $(z_{1}, z_{2})\mapsto (\eta z_{1}, \eta^{-1} z_{2})$.  Let
$H\subset \bZ/12\times \bZ/12$ be the subgroup of pairs $(a,b)$ with $a=b
\ {\rm mod}\  2$ and take the obvious action of $H$ on $\bC^{2}\times\bC^{2}$.
We consider the space of polynomials of bidegree (6,6) invariant under $H$.
 A basis for this space, written  in our usual  affine coordinates, is  
$ 1, s^{6}, t^{6},  s^{6}t^{6}, s^{3}t^{3} $
and there is an action of a 2-torus on the space, generated by $s\rightarrow
\mu s, t\rightarrow \nu t$.
We want to consider the corresponding curves in $\bP^{1}\times \bP^{1}$.
There is a similar discussion to the case of surfaces in subsection 6.1:
if any of the the coefficients of the first four monomials vanish we get
an unstable point for the torus action,  which we omit. Then we can use the
torus action to put the equation in the form
\begin{equation} P(s^{6}+t^{6}) + Q(1+s^{6}t^{6}) + R s^{3}t^{3}=0. \end{equation}
So we have a family of smooth curves $C_{P,Q,R}$ parametrised by an open
subset of $\bC\bP^{2}$. There is some residual equivalence: taking $s$ to
$-s$ gives $C_{P,Q,R}=C_{P,Q,-R}$ and taking $s$ to $s^{-1}$ gives $C_{P,Q,R}=C_{Q,P,R}$.
We write
$$   \cW^{H}= \bC\bP^{2}/(\bZ/2\times\bZ/2)$$ for the quotient by this action
of $\bZ/2\times\bZ/2$ and then the  moduli space of these smooth curves is
a subset of $\cW^{H}$.

\

Now go back to our $G$-invariant sextic surfaces  $S_{AB}$. The group action
gives a distinguished pencil of curves in  the canonical system  $\cO(2)$,
 defined by
$x_{2}y_{2}- \lambda  x_{1} y_{1}$. For $\lambda\neq 0,\infty$ the curve
is the intersection of $S_{AB}$ with a smooth quadric which we parametrise
by $x_{1}=s, y_{1}=t, x_{2}= st, y_{2}= \lambda$. Then the curve has equation
$$  s^{6}+t^{6}+ s^{6}t^{6}+ \lambda^{6}+ f_{AB}(\lambda) s^{3}t^{3}=0, $$
where $$f_{AB}(\lambda)= (1+\lambda)\left( A(1+\lambda)^{2}+ B(1-\lambda)^{2}\right).$$
Replacing $s, t$ by $\lambda^{1/2} s, \lambda^{1/2} t$ we put this curve
into our standard form
$$ \lambda^{3}(s^{6}+t^{6}) + (1+s^{6}t^{6}) + f_{AB}(\lambda) s^{3}t^{3}=0,$$

so $P=\lambda^{3}, Q=1, R=f_{AB}(\lambda)$. We view this as a degree $3$
map
$\Gamma_{A,B}: \bC\bP^{1}\rightarrow \bC\bP^{2}$. Let $\tau$ be the involution
of $\bC\bP^{2}$ defined by interchanging $P,Q$. It is easy to check that
the $\Gamma_{AB}$ are exactly the degree $3$ maps $\Gamma$ such that
\begin{enumerate} \item $\Gamma(\lambda^{-1})= \tau\Gamma(\lambda)$;
  \item $\Gamma(0)$ lies in the line $\{ Q=0\}\subset \bC\bP^{2}$ and $\Gamma$
has second order contact with the line at that point  (i.e. $Q\circ \Gamma=O(\lambda^{3})$
as $\lambda\rightarrow 0$); 
\item $\Gamma(\infty)$ lies in the line $\{ P=0\}\subset \bC\bP^{2}$ and
$\Gamma$
has second order contact with the line at that point.
\item $\Gamma$ does not pass through the point $P=Q=0$.
\item $\Gamma(-1)$ lies in the line $\{R=0\}$
\end{enumerate}
Of course the third item is a consequence of the first two.

Thus we see that, roughly speaking, our moduli space of sextic surfaces $S_{AB}$
can be interpreted as a  space of maps $\Gamma$ to a moduli space of curves
and we would like to investigate how this interacts with moduli space compactifications.

\

If $P$ and $Q$ are non-zero the curve defined by the equation (22) has at
worst ordinary double points so the question is how to extend the family
when $P$ or $Q$ vanish. We use the same procedure as before. For parameters
$\alpha,\beta$
we consider the complete intersection $C^{\alpha,\beta}$ defined by the equations
$$  x_{1}y_{1}= \alpha h\ \  , \  \ x_{2}y_{2}=\beta h\ \  ,\ \  x_{1}^{6}+y_{1}^{6}+
x_{2}^{6}+ y_{2}^{6}+ h^{3}=0$$ 
in $\bP^{4}_{w}$. If $\alpha,\beta$ are both nonzero then we view this as
a curve in the smooth quadric $x_{1}y_{1}= (\alpha/\beta) x_{2}y_{2}$ and
one finds that it is equivalent to $C_{P,Q,R}$ with $R=1, P=\alpha^{3} Q=\beta^{3}$.
When $\alpha=\beta=0$ we get a curve $C_{0}$ in the cone $\Cone(\Lambda)$
over the line configuration $\Lambda$. Let $\Sigma$ be the $\bZ/3$-cyclic
cover of $\bP^{1}$ branched over the six roots of  $z^{6}+1=0$, so $\Sigma$
has genus $2$ and there are three points in $\Sigma$ lying over $z=0$ and
three lying over $z=\infty$. The component of $C_{0}$ in each component of
$\Cone(\Lambda)$ is a copy of $\Sigma$ so we can obtain $C_{0}$ by taking
$4$ copies of $\Sigma$ and identifying  24 points (6 in each copy) in pairs.
So $C_{0}$ is a  stable curve in the sense of Deligne-Mumford, with $12$
ordinary double points. When just one of $\alpha, \beta$ is zero we are considering
the degeneration of the quadric to a pair of planes and we get a $1$-parameter
family of stable curves which deforms $C_{0}$ by smoothing $6$ of the double
points. 
 
 The conclusion is that in this case the naive compactification is the right
thing to consider. For each point $(P,Q,R)$ in $\bC\bP^{2}$ we have a stable
curve $C_{P,Q,R}$
and we obtain a moduli space $\cW^{H}= \bC\bP^{2}/(\bZ/2\times\bZ/2)$ of
stable
curves with this symmetry group,  contained in the full Deligne-Mumford compactification
$\ocM_{25}$ of curves of genus $25$. Our moduli space $M= \bC^{2}/\pm 1$
of surfaces $S_{AB}$ parametrises a family of curves $\uGamma_{AB}:\bC\bP^{1}\mapsto
\cW^{H}$. That is, we use the equivariance property (1) above to factor
$$ \bC\bP^{1}\rightarrow \bC\bP^{2}\rightarrow \cW^{H}$$
through the quotient $\bC\bP^{1}\rightarrow\bC\bP^{1}$ mapping $\lambda$
to $\lambda+\lambda^{-1}$, which identifies $\lambda$ with $ \lambda^{-1}$.
But it is easier to compute with the lifted maps $\Gamma_{AB}$.

On thing which is clarified by our family $C^{\alpha,\beta}$ is the second
order contact condition (2),(3) above. The curves $ C^{0,\beta}$ have an
additional $\bZ/3$-symmetry which means that this contact condition is the
condition that there is genuine family of curves $\Gamma(\lambda)$ for small
$\lambda$: this is just the fact that $P=\alpha^{3}$. We should really view
$\bC\bP^{2}$ as an orbifold with orbifold model at the origin given by the
quotient
$ \bC^{2}/(\bZ/3\times \bZ/3)$ 
 and this orbifold structure encodes the contact conditions (2),(3).

There is now another compactification $\oM_{\maps}$ of $M$ defined by the
theory of stable maps $\uGamma:\bC\bP^{1}\rightarrow \cW^{H}$. We will not
attempt to work this out in full here but one can observe some phenomena.
\begin{itemize}
\item For a surface on the boundary component from $III$ to $IV$ in $\oM$
defined by an equation 
$$  h^{3}+ c s^{2}t^{2} h +  (1+ s^{6}+t^{6}+s^{6} t^{6})=0$$
we take the pencil $st-\lambda h$ which gives the curves
$$ (\lambda^{-3}+ c\lambda^{-1})s^{3}t^{3} +(1+s^{6}+t^{6}+ s^{6|}t^{6})=0$$
in other words 
$$ P= Q= \frac{\lambda^{3}}{ 1+c \lambda^{2}} \ \ \ R=1; $$
a degree $3$ map to the line $P=Q$ in $\bC\bP^{2}$.
\item For a surface on the boundary component from $II$ to $IV$ in $\oM$
defined by an equation
$$  st h^{2}+ c s^{3}t^{3} + (1+ s^{6}+t^{6}+s^{6} t^{6})=0$$
we get
$$ P=Q= \frac{\lambda^{2}}{1+c\lambda^{2}}, \ \ \ R=1; $$
a degree 2 map to the line $P=Q$. 
\item For the surface $S_{IV}$ we have a pencil defined by $h^{+}-\lambda
h^{-}$ but for all but a finite number of $\lambda$ the curves are isomorphic
to $C_{0}$. 
\item For finite $A,B$ with $A+B\neq 0$ the image of the map $\Gamma_{AB}$
meets the line  $R=0$ at the points $P/Q= -1, \lambda_{+}^{3}, \lambda_{-}^{3}$
where
$$ \lambda_{\pm}= \frac{(B-A)\pm \sqrt{-4AB}}{A+B}. $$
This expression is homogeneous in $A,B$, so the $1$-parameter family of maps
$\Gamma_{tA,tB}$ meet the line at infinity in the same three points. These
points are recorded in the stable maps limit of the $\Gamma_{tA,tB}$ as $t\rightarrow
\infty$.
\end{itemize}
The first two items suggest that the compactification $\oM_{\maps}$ is not
the same as $\oM_{GIT}$ and the last shows that it is not the same as  $\oM$.
It seems likely that $\oM_{\maps}$ is a toric blow-up of each of these, with
 collapsing maps
$$  \oM_{GIT}\leftarrow \oM_{\maps}\rightarrow \oM . $$
  
Of course the fact that $\oM_{\maps}$ is different from $\oM$ does not rule
out the possibility of relating our enumerative theory to curve-counting
theories.

\

Counting parameters shows that the virtual (complex) dimension of the space
of maps $\Gamma$ satisfying the constraints (1)-(5) is $2$, the same as the
actual dimension, and so not the same as the virtual dimension of $\oM$.
 The explanation for this is similar to what we saw in Section 4, Remark
(7). The subset $\cW^{H}\subset \ocM_{25}$ lies
inside a larger family
$\ocM_{25}^{H}\subset \ocM_{25} $ of curves with an $H$-action.  More generally
there is a  moduli space $\cW$ of curves in the linear system $\cO(6,6)$
on $\bP^{1}\times \bP^{1}$ divided by $PSL(2,\bC)\times PSL(2,\bC)$ which
has dimension $49-1-6=42$, while the full moduli space $\cM_{25}$  of curves
of genus $25$ has dimension $3.25-3= 72$.  So for such a curve $C$ there
is a $30$-dimensional family of deformations which do not embed in $\bP^{1}\times\bP^{1}$.
To see this explicitly, we take the exact sequence of bundles on $C$
$$ 0\rightarrow TC\rightarrow \cO(2,0)\oplus\cO(0,2) \rightarrow \cO(6,6)\rightarrow
0$$
which gives
$$  H^{0}(C;\cO(6,6))\rightarrow H^{1}(TC)\rightarrow H^{1}(C, \cO(2,0)\oplus
\cO(0,2))\rightarrow H^{1}(C;\cO(6,6)). $$
The last term vanishes and the first term represents the deformations of
$C$ within $\bP^{1}\times \bP^{1}$. In terms of our $2$-dimensional spaces
$V^{\pm}$, as considered in subsection 6.3,  one finds that
$$  H^{1}(C; \cO(2,0)\oplus \cO(0,2))= s^{4}(V_{+})\otimes s^{2}(V_{-})\oplus
s^{2}(V_{-})\otimes s^{4}(V_{+})$$
which indeed has dimension $5.3+3.5=30$.

In our context we want to consider an $H$-invariant curve $C$ and the $H$-invariant
subspace of $s^{4}(V_{+})\otimes s^{2}(V_{-})\oplus  s^{4}(V_{+})\otimes
s^{2}(V_{-})$. It is easy to check that this space is $2$-dimensional so
$\cW^{H}$ has codimension-$2$ in  $\ocM_{25}^{H}$. To relate the deformation
theories it is clearest to work in the general case of $\cW\subset \ocM_{25}$
and a curve $\Gamma:\bP^{1}\rightarrow \cW$ defined by intersecting a sextic
surface $S$ with a pencil of quadrics in $\bP^{3}$. Without loss of generality
we  consider a point $p_{0}$ of $\bP^{1}$ corresponding to our standard quadric
$\bP^{1}\times \bP^{1}$ and a non-zero  tangent vector $v\in T\bP^{1}_{p}$.
The pencil involves another quadric and the choice of $v$ gives an element
 $$E_{v}\in H^{0}(\bP^{1}\times \bP^{1}; \cO(2,2))= s^{2}(V_{+}^{*})\otimes
s^{2}(V_{-}^{*}). $$ Now suppose that we have an element $\Omega\in \Lambda^{2}U^{*}$.
As we recalled in subsection 5.3 we can regard $\Omega$ as lying in $s^{2}(V^{*}_{+})\oplus
s^{2}(V^{*}_{-})$.
Then we have a product 
$$  \Omega. E_{v}\in s^{4}(V^{*}_{+})\otimes s^{2}(V^{*}_{-})\oplus
s^{2}(V^{*}_{-})\otimes s^{4}(V^{*}_{+}), $$
which is the dual of $H^{1}(C; \cO(2,0)\oplus \cO(0,2))$. Now, as we have
explained above,  this last space can be regarded as the fibre of the normal
bundle $N_{\cW}$ of $\cW$ in $\ocM_{25}$. The upshot is that we get map 
$$\Lambda^{2}U^{*}\rightarrow H^{0}(\bP^{1}, T^{*}\bP^{1}\otimes \Gamma^{*}(N_{\cW}^{*})).$$
Using Serre duality on $\bP^{1}$ the transpose is a map
$$  
 H^{1}(\bP^{1}, \Gamma^{*}(N_{\cW}))\rightarrow \Lambda^{2} U. $$
 Composing with the map induced by projection of $T\ocM_{25}$ to the normal
bundle we get a map from $H^{1}(\bP^{1}, \Gamma^{*}T\ocM_{25})$ to $\Lambda^{2}
U$ which relates the obstruction spaces in the two theories.  

\

As we wrote at the beginning of this section,  the main motivation for studying
this example with finite group action is to gain insight into the larger
questions, such as for general sextic surfaces. The dimension of the space
of pairs consisting of a sextic surface $S$ and a pencil in $\vert K_{S}\vert$
is $68+ 16=84$ and for each such pair we get (roughly speaking) a rational
curve $\Gamma$ in the 42-dimensional space $\cW$ inside the $72$-dimensional
$\ocM_{25}$.  Since $84=2.42$ we expect that for typical points $C_{1}, C_{2}\in
\cW$ there is a $1$-dimensional space of curves $\Gamma$ through $C_{1}$
and $ C_{2}$. But if we consider the space $\cW$ as a subset of $\ocM_{25}$
we expect that the virtual dimensions reduce by $6$. It seems interesting
to study both the detailed geometry and the enumerative geometry of this
situation.  

{\bf Acknowledgements}  The author is very grateful to Richard Thomas for
many suggestions and comments on a preliminary version of this article and
to Ivan Smith for  pertinent remarks some years ago. The author's work is supported by the Simons Foundation through the
Simons Collaboration on Special Holonomy.



\begin{thebibliography}{99}
\bibitem{kn:Alexeev} V. Alexeev {\em Moduli spaces $M_{g,n}(W)$ for surfaces}
in : Higher dimensional complex varieties Trento 1994 de Gruyter (1996) 1-22
\bibitem{kn:Atiyah2} M. Atiyah {\em Riemann surfaces and spin structures}
Ann. Sci. Ecole Norm. Sup. 4 (1971) 47-62
\bibitem{kn:AB} M. Atiyah and R. Bott {\em the Yang-Mills equations on Riemann
surfaces} Phil. Trans. Roy. Soc. London A 308 (1982) 523-615
\bibitem{kn:AHS} M. Atiyah, N. Hitchin and I. Singer {\em Self-duality in
four-dimensional Riemannian geometry} Proc. Roy. Soc. London A 362 (1978)
425-61
\bibitem{kn:BF} S. Bauer and M. Furuta {\em A stable cohomotopy refinement
of Seiberg-Witten invariants,I} Inventiones Math. 155 (2004) 1-19
\bibitem{kn:BPV} W. Barth, C. Peters and A. Van de Ven  {\em Compact complex
surfaces} Springer 1984
\bibitem{kn:BFan} K. Behrend and B. Fantechi {\em The intrinsic normal cone}
Inventiones Math. (1997) 45-88
\bibitem{kn:Berman} R. Berman and H. Guenancia {\em K\"ahler-Einstein metrics
on singular varieties and log canonical pairs} Geometric and Functional Analysis
24 (2014) 1683-1730
\bibitem{kn:BJ} D. Borisov and D. Joyce {\em Virtual fundamental classes
for moduli spaces of sheaves on Calabi-Yau fourfolds} Geometry and Topology
21 (2017) 3231-3311
\bibitem{kn:Cao} Y. Cao and N. Leung {\em Donaldson-Thomas theory for Calabi-Yau
4-folds } arxiv 1407.7659
\bibitem{kn:Catanese} F. Catanese {\em On the moduli spaces of surfaces of
general type} Jour. Differential Geometry 19 (1984) 483-515
\bibitem{kn:Donaldson1} S. Donaldson {\em Polynomial Invariants of smooth
four-manifolds} Topology 29 (1990) 257-315
\bibitem{kn:DF} S. Donaldson and R. Friedman {\em Connected sums of self-dual
manifolds and deformations of singular spaces} Nonlinearity 2 (1989) 197-239
\bibitem{kn:DT} S. Donaldson and R. Thomas {\em Gauge theory in higher dimensions}
In: The Geometric Universe S. Huggett {\it et al} eds. Oxford U.P. (1998)
31-47
\bibitem{kn:FL} P. Feehan and T. Leness {\em Witten's conjecture for many
4-manifolds of simple type} Jour. European Math. Soc. 17 (2015) 899-923
\bibitem{kn:FO} M. Furuta and H. Ohta {\em Differentiable structures on punctured
4-manifolds} Topology and Applications 51 (1993) 291-301
\bibitem{kn:GRW} S. Galatius and O. Randall-Williams {\em Moduli spaces of
manifolds: a user's guide} in:Handbook of homotopy theory, Routledge 2019
arxiv 1811.08151
\bibitem{kn:Hacking} P. Hacking {\em Compact moduli spaces of surfaces of
general type} arxiv 1107.2717
\bibitem{kn:Vetal} H-J. Hein, S. Sun, J. Viaclovsky and R. Zhang {\em Nilpotent
structures and collapsing Ricci-flat metrics on K3 surfaces} arxiv 1807.09367
\bibitem{kn:JK} L. Jeffrey and F. Kirwan {\em Localisation for nonabelian
group actions} Topology 34 (1995) 291-327
\bibitem{kn:KK} A. King and D. Kotschick {\em The deformation theory of ant-self-dual
conformal  structures} Math. Annalen 294 (19920 591-610
\bibitem{kn:Kobayashi} R. Kobayashi {\em Einstein-K\"ahler V-manifolds on
open Satake V-surfaces with isolated quotient singularities} Math. Annalen.
272 (1985) 385-98
\bibitem{kn:KSB} J. Koll\'ar and N. Shepherd-Barron {\em Threefolds and deformations
of surface singularities} Inventiones Math. 91 (1988) 299-338
\bibitem{kn:Laza} R. Laza {\em The KSBA compactification of the moduli space
of degree 2 K3 pairs} Jour. European Math. Soc. 225-79
\bibitem{kn:LT} J. Li and G. Tian {\em Virtual moduli cycles and Gromov-Witten
invariants of algebraic varieties}  Jour Amer. Math. Soc. 11 (1998) 119-174
\bibitem{kn:LiLiu} T.J. Li and A-K. Liu {\em Family Seiberg-Witten invariants
 and wall crossing formulas} Commun. Anal. Geom. 9 (2001) 777-823
\bibitem{kn:Manetti} M. Manetti {\em Iterated double covers and conected
components of moduli spaces of general type} Topology 36 (1997) 745-64
\bibitem{kn:Nak} N. Nakamura {\em The Seiberg-Witten equations for families
and diffeomorphisms of 4-manifolds} Asian Jour. Math. 7 (2003) 133-138
\bibitem{kn:PXu} Z. Patakfalvi and C. Xu {\em Ampleness of the CM line bundle
on the moduli spaces of canonically polarised varieties} Algebraic Geometry
4 (2017) 29-39
\bibitem{kn:Poon} Y. Sun Poon {\em Compact self-dual manifolds with positive
scalar curvature} Jour. Differential Geometry 24 (1986) 97-132
\bibitem{kn:RW} O. Randall-Williams {\em Miller-Mumford-Morita classes vanish
on moduli spaces of hypersurfaces} available at https://www.dpmms.cam.ac.uk/~or257/
\bibitem{kn:Smale} S. Smale {\em An infinite-dimensional version of Sard's
Theorem} Amer. Jour. Math. 87 (1965) 861-66
\bibitem{kn:Taubes1} C. Taubes  {\em Casson's invariant and gauge theory}
Jour. Differential Geometry  31 (1990) 547-599
\bibitem{kn:Taubes2} C. Taubes {\em The behaviour of sequences of solutions
to the Vafa-Witten equations} arxiv 1702.04610
\bibitem{kn:Thomas}  R. Thomas {\em A holomorphic Casson invariant for Calabi-Yau
3-folds and bundles on K3 fibrations} Jour. Differential Geometry 54 (2000)
367-438
\bibitem{kn:TV} G. Tian and J. Viaclovsky{\em Moduli spaces of critical Riemannian
metrics in dimension four} Advances in Math. 196(2005) 346-372
\end{thebibliography}
\end{document}